\documentclass[11pt]{article}
\usepackage[T1]{fontenc}
\usepackage{latexsym,amssymb,amsmath,amsfonts,amsthm,array}
\usepackage{graphics,graphicx,mathrsfs}
\usepackage{color}
\usepackage{float} 
\usepackage{subfigure} 
\usepackage{url,epsfig,tikz}
\usepackage{pgfplots}
\usepackage{epstopdf}
\usetikzlibrary{shapes.geometric}

\newcommand{\R}{{\mat R}}

\newcommand{\N}{{\mat N}}
\newcommand{\C}{{\mat C}}
\newcommand{\Sb}{{\mat S}}

\newcommand{\be}{\begin{eqnarray}}
\newcommand{\ben}{\begin{eqnarray*}}
\newcommand{\en}{\end{eqnarray}}
\newcommand{\enn}{\end{eqnarray*}}
\newcommand{\ba}{\begin{aligned}}
\newcommand{\an}{\end{aligned}}

\newcommand{\wt}{\widetilde}

\newcommand{\B}{{\rm B}}

\newcommand{\mat}{\mathbb}

\newcommand{\supp}{{\rm supp}\,\,}

\newtheorem{thm}{Theorem}[section]

\newtheorem{lem}{Lemma}[section]

\newtheorem{assum}{Assumption}[section]

\definecolor{rot}{rgb}{1,0,0}
\definecolor{hw}{rgb}{0,0,1}

\begin{document}
\renewcommand{\theequation}{\arabic{section}.\arabic{equation}}
\title{\bf Uniqueness, stability and algorithm for an inverse wave-number-dependent source problem}
\author{
Mengjie Zhao\thanks{School of Mathematical Sciences and LPMC, Nankai University,
300071 Tianjing, China. ({\tt 1120200028@mail.nankai.edu.cn})}
\and
Suliang Si\thanks{(Corresponding author) School of Mathematics and Statistics, Shandong University of Technology,
255049 Shandong, China. ({\tt sisuliang@amss.ac.cn})}
\and
Guanghui Hu\thanks{School of Mathematical Sciences and LPMC, Nankai University,
300071 Tianjing, China. ({\tt ghhu@nankai.edu.cn})}
}

\date{}


\maketitle

\begin{abstract}
This paper is concerned with an inverse wavenumber/frequency-dependent source problem for the Helmholtz equation. In two and three dimensions, the unknown source term is supposed to be compactly supported in spatial variables but independent on one spatial variable. The dependence of the source function on wavenumber/frequency is supposed to be unknown. Based on the Dirichlet-Laplacian and Fourier-Transform methods, we develop two efficient non-iterative numerical algorithms to recover the wavenumber-dependent source.
Uniqueness proof and increasing stability analysis are carried out in terms of the boundary measurement data of Dirichlet kind. Numerical experiments are conducted to illustrate the effectiveness and efficiency of the 
proposed methods.

 
\end{abstract}

\vspace{.2in} {\bf Keywords: inverse source problem, wavenumber-dependent source, Helmholtz equation, increasing stability.
}

%

\section{Introduction and motivations}
\subsection{Statement of the problem}

\par Consider a time-dependent acoustic wave radiating problem in $\R^d$ ($d=2,3$) with the homogeneous initial conditions 
\be\label{equ:ut1}
\left\{
	\ba
		\partial^2_t U(x,t)-\Delta U(x,t)= - S(x,t),  \quad  x&\in \R^d,\,t>0,\\
		U(x,0)=\partial_{t} U(x,0)=0, \quad  x&\in \R^d ,
	\an
	\right. 
\en
where $x=(\wt x, x_d)$, $\wt x=(x_1, x_2,\cdots, x_{d-1})\in \R^{d-1}$ and  $U(x,t)$ represents the wave fields. To state the motivation of our problem,  we suppose that the source term $S(x,t)$ is excited by a moving point source along the trajectory $a(t):\R_+\longmapsto \R^d$ for $t\in [0,T]$. More precisely, we suppose that the source function takes the form (see e.g., \cite{HK2020})
\be\label{equ:source} S(x,t)=\delta(x-a(t))\chi (t),\en
where the symbol $\delta$ is the Dirac delta function and $\chi (t)$  is the characteristic function given by  
\ben
	\chi (t)=\left\{
	\begin{array}{lll}	
		1, \qquad t \in {[ 0,T]  },\\ 
		0, \qquad t \notin {[0,T]}.
	\end{array}
	\right. 	
\enn
Assume further that the point source moves on the hyperplane $x_d=h$ for some $h>0$, so that the trajectory function can be written as $a(t)=(\wt a(t),h)$, where $\wt a(t):\R_+\longmapsto \R^{d-1}$ is supposed to be a smooth function.  
Then the source function can be written as $S(x,t)=\delta(\wt x-\wt a(t))\,\chi(t)\,\delta(x_d-h)$. Physically, the Delta function can be approximated by functions with a compact support in the distributional sense. Hence the source term $S(x,t)$
can be physically approximated by 
$F(\tilde{x},t) g(x_d)$,
where $F$ and $g$ have compact supports. 
Below we shall identify $S(x,t)$ with $F(\tilde{x},t)g(x_d)$ for the convenience of mathematical analysis.
The solution of problem \eqref{equ:ut1} can be given explicitly as 
\ben\label{equ:solutionUt}
U(x,t)=-\int_{\R^{+}}\int_{\R^{d-1}}G_d(x-y;t-\tau)S(y,\tau)dyd\tau,\quad  x\in \R^d,\, t>0,
\enn
where 
\ben
G_d(x;t)=\left\{
\ba
&\dfrac{H(t-|x|)}{2\pi \sqrt{t^2-|x|^2}}, &\qquad  d=2, \quad |x|\neq t, \\
&\dfrac{\delta(t-|x|)}{4\pi |x|},         &\qquad  d=3, \quad |x|\neq t,  \  x\neq 0,
\an
\qquad 
x\in \R^d \right.
\enn 
is the fundamental solution to the wave equation in $\R^d$ and $H(\tau)$ is the Heaviside function
\ben
	H(\tau)=\left\{
	\begin{array}{l}	
		1, \qquad \tau\geq 0,\\
		0, \qquad \tau<0.                 
		
	\end{array}
	\right. 	
\enn 
The $n$-dimensional ($n=1,2$) Fourier and inverse Fourier transforms in this paper are defined respectively by
\ben
(\mathcal{F}_\xi\,\rho)(\tau)&=&\dfrac{1}{(2\pi)^{n/2}}\int_{\R^n} \rho(\xi)e^{-i\tau\cdot\xi}\,d\xi,\quad \tau\in \R^n,\\
(\mathcal{F}_\tau^{-1} q)(\xi)&=&\dfrac{1}{(2\pi)^{n/2}}\int_{\R^n} q(\tau)e^{i\tau\cdot\xi}\,d\tau,\quad \,\,\,\xi\in \R^n,
\enn 
where the subscript denotes the variable to be Fourier or inverse Fourier transformed. 
The inverse Fourier transform of the source term takes the form 
\ben\label{equ:fft}
(\mathcal{F}_t^{-1}F)(\wt x,k)g(x_d):=f(\wt x,k)g(x_d).
\enn
It is obvious that the spatial function $\wt x\longmapsto f(\wt x,k)$ is compactly supported for all $k>0$. The temporal inverse Fourier transform of $U(x,t)$ is given by 
\be \label{equ:DP}
u(x,k):=(\mathcal{F}^{-1}_t U)(x,k)= - \dfrac{1}{\sqrt{2\pi}} \int_{\R^d}  \Phi_d(x-y;k)  f(\wt y,k)g(y_d) dy , \qquad x\in \R^d, \, k>0,
\en 
where $\Phi_d(x;k)$ is the fundamental solution  to the Helmholtz equation $(\Delta+k^2)u=0$ in $\R^d$, given by
\ben 
\Phi_d(x-y;k)=(\mathcal{F}^{-1}_t G_d)(x-y;k)=\left\{
\ba
&\dfrac{i}{4}H^{(1)}_0(k|x-y|),   &\quad d=2,\\
&\dfrac{e^{ik|x-y|}}{4\pi |x-y|}, \   &\quad d=3,
\an
\qquad 
x,y\in \R^d, \  x\neq y.\right.
\enn 
Here $H_0^{(1)}$ is the first kind Hankel function of order zero. Applying the inverse Fourier transform w.r.t. the time variable to the wave equation \eqref{equ:ut1},  we conclude that
\be\label{equ:uk}
\Delta  u(x,k)+k^{2} u(x,k) = f(\wt x,k)g(x_d), \qquad x\in \R^d, \,  k>0.
\en
Moreover, the radiated field $u(x,k)$ satisfies the following Sommerfeld radiation condition
\be\label{equ:SRC}
\lim\limits_{r\rightarrow\infty} r^{\frac{d-1}{2}}(\partial_r u-iku)=0,   \quad  r=|x|,
\en
which holds uniformly in all directions $\hat{x}=x/{|x|}\in \Sb^{d-1}:=\{x\in \R^d:|x|=1\}$.

\par Denote by $B_R\subset \R^d$ the sphere centered at the origin with the radius $R>0$ and by $\partial B_R:=\{x\in \R^d:\ |x|=R \}$ the boundary of $B_R$.
The unit normal vector $\nu$ at the boundary $\partial B_R$ is supposed to direct into the exterior of $B_R$. 
 
We make the following assumptions on the source functions $f(\wt x,k)$ and $g(x_d)$:
\begin{assum}\label{assumA} It holds that 
	\begin{itemize}
		
		\item[$(a)$] There exist non-empty sets $\wt D\subset \R^{d-1}$ and $H\subset \R_+$ such that $\supp g=\overline{H}$, $\supp f(\cdot,k)=\overline{\wt D}$, $(\wt D\times H)\subset\subset B_R$.
		
		\item[$(b)$] $f(\cdot ,k)\in L^2(\wt B_R)$, where $\wt B_R:=\{\wt x\in \R^{d-1}:\ |\wt x|<R\}$ and $g\in L^2(-R,R)$.
		
		\item[$(c)$] $f(\wt x,k)$ is analytic with respect to $k$,  $g$ is a priori known and $g\not\equiv 0$.

	\end{itemize}   
\end{assum}

In this paper, we are interested in the following inverse problem:

\noindent\textbf{(IP)}:\ Determine the source function $f(\wt x,k)\in L^2(\wt B_R)$ for all $k>0$ from knowledge of the multi-frequency near-field measurement \be\label{equ:data}
\big\{ u(x,k):\ x\in \partial B_R, \ k\in [k_{\min},k_{\max}]\},
\en
where $[k_{\min},k_{\max}]$ denotes an interval of available wavenumbers with $0<k_{\min}<k_{\max}$.

\subsection{Literature review}
Inverse source problems in wave propagation are of great importance in various scientific and
engineering fields such as antenna design
and synthesis, medical imaging and photo-acoustic tomography \cite{A1999, FKM2004, H2000, SU2011}.
Consequently, a great deal of mathematical and numerical results are available, especially for wavenumber-independent source problems of
the Helmholtz equation. In general,
there is an obstruction to uniqueness for inverse source problems with single-frequency data due to the existence of non-radiating sources  \cite{BC1977, I2017}. Computationally, a more serious issue is the stability analysis, i.e., to quantify how does a small variation of the data lead to a huge error in the reconstruction. Hence it is crucial to study the stability of inverse source problems. Recently, it has been realized that the use of multi-frequency data is an effective approach to overcome the difficulties of non-uniqueness and instability which are encountered at a single frequency.   In \cite{BLT2010}, Bao et al.
initialized the mathematical study on the stability of inverse wavenumber-independent source problems for the Helmholtz equation
by using multi-frequency data. In the recent works \cite{BLL2015, BLT2010, CIL2016, EV2009, HKZ2020}, uniqueness and stability results have been proved for recovering source terms from multi-frequency boundary measurements. In \cite{BLT2010}, the authors treated an interior inverse source problem for the Helmholtz equation from boundary Cauchy data for multiple wave numbers and showed an increasing stability result. 
A uniqueness result and numerical algorithm for recovering the location and shape of a supported acoustic source function from boundary measurements at many frequencies were shown in \cite{CIL2016}. The approach of \cite{EV2009} also provides inspirations for establishing uniqueness result. It was shown that an acoustic source of the form $S(x,k)=f(x)g(k)$ can be identified uniquely by boundary measurements of the acoustic field. The proof relies on complete sets of solutions to the homogeneous Helmholtz equation. An error estimate and a numerical algorithm for the solution are presented. Several numerical methods for solving the multi-frequency inverse source problem have been proposed. We refer the reader to \cite{BLRX2015, ZG2020} for the Fourier method and recursive algorithm. In addition, there are other uniqueness, increasing stability and algorithm results; see also \cite{ACTV2012,AHLS2020,BLT2021,BB2017,CH2018,E2018,EI2020,  IL2018,I2007, INU2014,LY2017,LZZ2020} and the references therein.

To the best of our knowledge, there are only few works on  inverse wave-number-dependent source problems. Even using multi-frequency data, there is no uniqueness in recovering general wave-number-dependent (or frequency-dependent) source functions. Difficulties arise from the fact that the far-field measurement data are no longer the Fourier transform form of the source function. The wave-number-dependent source functions are closely connected to time-dependent source functions in the time domain \cite{GH24,HK2020, HKZ2020}. We refer to \cite{GH24,GHZ23} for sampling methods to image the support of a special class of wave-number-dependent source functions with known radiating periods.   In this paper, we consider acoustic source functions of the form $f(\wt x,k)g(x_d)$ in the frequency regime. The unknown wave-number-dependent source function $f$ is independent of one spatial variable and the source term $g(x_d)$, $x_d\in\R$ is supposed to be a priori known. Motivated by existing results for wave-number-independent sources \cite{BLT2010,CIL2016,EV2009,I2007,LY2017}, we prove uniqueness and increasing stability in covering  
 the source function $\wt x\mapsto f(\wt x,k)$, $\wt x\in\R^{d-1}$ from the Dirichlet boundary data \eqref{equ:data} at a fixed frequency/wavenumber $k>0$. This implies that the multi-frequency/wavenumber data can be used to uniquely identify the function $f$ by analytic continuation.
Inspired by the idea of \cite{EV2009}, we choose a complete orthonormal set to prove that the source function $f$ can be identified uniquely by the Dirichlet data measured on the boundary; see the uniqueness result shown in Theorem \ref{thm:DL}.
The uniqueness proof and inversion algorithms are based on 
two methods: Fourier transform method and Dirichlet-Laplacian method. The increasing stability estimates are derived from  the Fourier transform method and the technique of \cite{CIL2016}.  

\par The remaining part of the paper is organized as follows. In Section 2, we show that the source term $f(\wt x,k)$ can be identified uniquely by the multi-frequency near-field data. Section 3  is devoted to the increasing stability analysis. Numerical experiments are presented in Section 4 to verify the effectiveness of the proposed method.

\section{Uniqueness}\label{sec:unique}
This section is concerned with uniqueness of the inverse wave-number-dependent source problem by measuring the near-field acoustic field at many frequencies. 
In subsection \ref{subsec:F} we apply the Fourier transform method with test functions, while subsection \ref{subsec:D} demonstrates the feasibility of achieving uniqueness throught the Dirichlet-Laplacian method. The outcomes of our investigation are presented in Theorems \ref{thm:FT} and \ref{thm:DL} respectively.

\par

\subsection{Uniqueness via Fourier transform}\label{subsec:F} 
\par In this subsection, we aim to prove uniqueness by applying the Fourier transform method, which is stated as follows.
\begin{thm}\label{thm:FT}
Suppose that the source functions $f$ and $g$ satisfy the Assumption \eqref{assumA} and that $u(\cdot,k)\in H^2(B_R)$ for all $k>0$ is the unique solution to the inhomogeneous Helmholtz equation \eqref{equ:uk} with the Sommerfeld radiation condition \eqref{equ:SRC}. Then the unknown source $f(\wt x,k)\in L^2(\wt B_R)$ for $k>0$ can be uniquely determined by the Dirichlet data $\{u(x,k): x\in \partial B_R,\ k\in [k_{\min},k_{\max}]\}$.  
\end{thm}
\begin{proof} Since the inverse source problem is linear and $g$ is a priori known, we shall prove $f\equiv 0$ by assuming 
 $u(x,k)=0$ on $\partial B_R$ for all $k\in [k_{\min},k_{\max}]$.
 Using the uniqueness of the exterior boundary value problem of the Helmholtz equation \cite[Theorem 9.10, Theorem 9.11]{W2000}, 
 one can conclude that $u(x,k)=0$ in $\R^d \textbackslash \overline{B}_R$ and thus $\partial_\nu u(x,k)=0$ on $\partial B_R$ for $k\in [k_{\min},k_{\max}]$. Define the wavenumber-dependent  test functions 
\be\label{equ:test1} 
	\psi(x,k)=e^{-i\xi\cdot x},\quad \xi=(\wt \xi,\xi_d) \in \R^d, \quad |\xi|=k. 
	\en	
	For all $k>0$, it is easy to verify $\psi(x,k)$ satisfies the Helmholtz equation $ \Delta \psi(x,k)+k^2\psi(x,k)=0$ in $\R^d$.
Multiplying $\psi(x,k)$ on both sides of \eqref{equ:uk} and integrating  over $B_R$, we obtain 
	\be\label{equ:Green1}
	\begin{aligned}
		\int_{B_R} f(\wt x,k)g(x_d)\psi(x,k)dx &=\int_{B_R} \left(  \Delta u(x,k)+k^2 u(x,k)\right)   \psi(x,k) dx,\\ &=\int_{\partial B_R}\left[ \dfrac{\partial u(x,k)}{\partial \nu(x)} \psi(x,k)-u(x,k)\dfrac{\partial \psi(x,k)}{\partial \nu(x)}\right] ds\left( x\right), \\ &=0.
	\end{aligned}
    \en
   	Hence, we have
   	\ben
   	\begin{aligned}
   		\left(\int_{\wt D} f(\wt x,k)e^{-i\wt \xi \cdot \wt x}d\wt x\right)\;\left( \int_{H} g(x_d)e^{-i\xi_dx_d}dx_d\right)=0,   \quad  k\in[k_{\min},k_{\max}].   
   	\end{aligned}
   \enn
 Below we fixed some  $k\in [k_{\min},k_{\max}]$. 
 Noting that $\supp g(x_d)=\overline{H}$ by our 
Assumption \eqref{assumA}, one can always find some $0<\eta<k$ such that $\mathcal{F} g\neq0$ in $(0,\eta)$.
   Therefore, for all $k\in[k_{\min},k_{\max}]$ we have 
	\be\label{equ:FTformula}
		\int_{\wt D} f(\wt x,k)e^{-i\wt \xi \cdot \wt x}d\wt x=0 \quad\mbox{for all}\quad \wt \xi \in S_k,
	\en
where  $S_k:=\{\wt x\in\R^{d-1}:\ |\wt x|^2=k^2-x_d^2,\ x_d\in(0,\eta) \}$.
From \eqref{equ:FTformula}, we deduce that $(\mathcal{F}_{\wt x}\,f)(\wt \xi,k)=0$ for all $\wt \xi\in S_k$ with  $k\in[k_{\min},k_{\max}]$. Hence, it holds that $(\mathcal{F}_{\wt x}\,f) (\wt \xi,k) =0$ for all $\wt \xi\in\R^{d-1}$, since $\mathcal{F}_{\wt x}\,f$ is analytic in $\R^{d-1}$ and $S_k$ is a non-empty open subset of $\R^{d-1}$. Taking the inverse Fourier transform on $(\mathcal{F}_{\wt x}\,f) (\wt \xi,k)$ yields $\mathcal{F}_{\wt \xi}^{-1}(\mathcal{F}_{\wt x}\,f)=f(\wt x,k)=0$ for all $\wt x \in\R^{d-1}$ with $k\in[k_{\min},k_{\max}]$. Therefore we conclude that $f(\wt x,k)=0$ with $\wt x\in\R^{d-1}$ for all $k>0$ by applying the analyticity of $f(\wt x,k)$ with respect to $k$.
\end{proof}

\subsection{Uniqueness via Dirichlet eigenfunctions}\label{subsec:D}

The proof presented in subsection \ref{subsec:F} cannot be directly applied to inhomogeneous background media, due to the difficulties in constructing corresponding test functions. Here we explore an alternative method based on the completeness of Dirichlet eigenfunctions for the negative Laplacian.
Let $v\in H_0^1(\wt B_R)$ be an eigenfunction of the Dirichlet Laplacian in $\wt B_R$, i.e., the solution to the Dirichlet problem for the homogeneous Helmholtz equation in $\wt B_R$, 
\ben
		-\Delta_{\wt x} v(\wt x) =\lambda^2 v(\wt x), \quad \wt x\in \wt B_R,\qquad
		v(\wt x) = 0, \, \quad \wt x\in \partial \wt B_R.
\enn
Here $\lambda^2$ is known as the Dirichlet eigenvalue for $\wt B_R$. Denote the sequence of eigenvalues as $\left\lbrace \lambda_{m}^2\right\rbrace _{m\in\mathbb{N}^+}$ counted with their multiplicity and the Dirichlet eigenfunctions $\lbrace v^{(m)}_j\rbrace _{j=1}^{m_j}$ with $m_j\in \mathbb{N}^+$   associated with $\lambda_{m}^2$. We renumber all eigenfunctions as $\lbrace v_n\rbrace_{|n|\in\mathbb{N}^+}$. Since the Dirichlet Laplacian operator is positive and self-adjoint in $L^2(\wt B_R)$, the set of eigenfunctions $\lbrace v_n\rbrace_{|n|\in\mathbb{N}^+}$ forms a complete orthonormal basis in $L^2(\wt B_R)$. Thus for all $h\in L^2(\wt B_R)$ we have the expansion
\be\label{equ:expansion}
h(y)=\sum\limits_{|n|\in\mathbb{N}^+}h_n\;  v_n(y), \quad y\in \wt B_R,
\en 
where the coefficients $h_n$ are given by
\ben\label{equ:coefficient}
h_n=\int_{\wt B_R} h(y)\,\overline{v}_n(y)\,d y, \quad |n|\in\mathbb{N}^+.
\enn


We remark that the only possible accumulating point of $\lambda_{n}^2$ is at infinity.
\begin{lem}\label{lem} Let $\left\lbrace \lambda_{n}^2\right\rbrace _{n=1}^{+\infty}$ be the sequence of Dirichlet eigenvalues over $\wt B_R$.
Then there exists a non-empty open subset $K\subset [k_{\min},k_{\max}]$, which satisfies $\lambda_n\notin K$ for all $n\in \N^{+}$, such that 
\ben
\int_{H}g(t)e^{\sqrt{\lambda_n^2-k^2}\,t}\,dt\neq0  \quad\mbox{for all}\quad k\in K
\enn
holds for all $n\in\mathbb{N}^+$.
\end{lem}

\begin{proof}
    Suppose that for any non-empty open subset $K\subset [k_{\min},k_{\max}]$, which satisfies $\lambda_n\notin K$ for all $n\in \N^{+}$, one can find an eigenvalue $\lambda_l$ with some $l\in \N^{+}$ such that 
    \be\label{equ:contradict}
    \int_{H}g(t)e^{\sqrt{\lambda_l^2-k^2}\,t}\,dt=0, \quad k\in K.
    \en
For the aforementioned $l\in \N^{+}$, we introduce two subsets of $K$:
    \ben 
       K_{1}:=\{k\in K:\ k^2<\lambda_l^2\}, \quad K_{2}:=\{k\in K:\ k^2>\lambda_l^2\}.
    \enn
    It is obvious that $K_{1}$ and $K_{2}$ cannot be both empty. We consider the following two cases.\\
    \textbf{Case 1}:
    \noindent  $K_{1}\subset K$ is non-empty. Introducing the variable $s:=\sqrt{\lambda_l^2-k^2}$ for $k\in K_1$, one obtains $s\in (c,d)$ for some $0\leq c<d<\infty$. Observing that the function $s\longmapsto \int_{H}g(t)e^{s\, t}\,dt$ is analytic with respect to $s\in\R$ and applying the Taylor expansion for the exponential function, we obtain
    \ben
      0=\int_{H}g(t)e^{s\,t}dt=\sum_{n\in\N}A_n s^n, 
      \qquad
      A_n:=\dfrac{1}{n!}\int_{H}g(t)\, {t}^n\,dt,&\quad n\in \N
    \enn
for all $s\in \R$, and hence $A_n=0$. This together with the completeness of polynomials in $L^2(H)$ implies  $g=0$, which   
is a contradiction to our assumption $g\not\equiv 0$.   
\\    
\textbf{Case 2}:
    \noindent  $K_{2}\subset K$ is non-empty. In this case one obtains

    \ben
    \int_{H}g(t)e^{\sqrt{\lambda_l^2-k^2}\,t}dt=\int_{H}g(t)e^{i\sqrt{k^2-\lambda_l^2}\,t}dt=(\mathcal{F}^{-1}_{t}g)(\sqrt{k^2-\lambda_l^2})=0, \quad  k\in K_{2}.
    \enn
Again using the analyticity of the inverse Fourier transform of $g(t)$, we conclude that $(\mathcal{F}^{-1}_{t}g)=0$ in $\R$. Therefore, we get $\mathcal{F}_k (\mathcal{F}^{-1}_{t}g)=g=0$ in $\R$, 
    which contradicts the assumption $g\not\equiv 0$. 
\end{proof}


The uniqueness proof for Theorem \ref{thm:DL} below differs from Theorem \ref{thm:FT} in the choice of test functions.
\begin{thm}\label{thm:DL}
Under the assumption of Theorem \ref{thm:FT},
the source function $f(\wt x,k)\in L^2(\wt B_R)$ for $k>0$ can be uniquely determined by the Dirichlet data $\{u(x,k): x\in \partial B_R,\ k\in [k_{\min},k_{\max}]\}$ .  
	
\end{thm}

\begin{proof}
Define the $k$-dependent test functions
\ben
\varphi_n(x,k)=v_{n}(\wt x)\;e^{\sqrt{\lambda_{n}^2-k^2}\,x_d},\quad  n\in\mathbb{N}^+,
\enn
where $\lambda_{n}^2$ are the Dirichlet eigenvalues and $v_n$ are the Dirichlet eigenfunctions associated with $\lambda_{n}^2$. For all $n\in\mathbb{N}^+$, the function $\varphi_n$ satisfies the Helmholtz equation $ \Delta \varphi_n(x,k)+k^2\varphi_n(x,k)=0$ in $\wt B_R$.

\par Assuming $u(x,k)=0$ on $\partial \B_R$ for $k\in [k_{\min},k_{\max}]$, we deduce that $\partial_\nu u(x,k)=0$ on $\partial \B_R$ for $k\in [k_{\min},k_{\max}]$.
Multiplying $\varphi_n(x,k)$ on both sides of \eqref{equ:uk} and integrating  over $B_R$ for all $n\in\mathbb{N}^+$, we obtain 
\be\label{equ:Green2}
\begin{aligned}
\int_{B_R} f(\wt x,k)g(x_d)\varphi_n(x,k)dx &=\int_{B_R} (\Delta u(x,k)+k^2 u(x,k))\varphi_n(x,k) dx\\ &=\int_{\partial B_R}\left(\dfrac{\partial u(x,k)}{\partial \nu(x)} \varphi_n(x,k)-u(x,k)\dfrac{\partial \varphi_n(x,k)}{\partial \nu(x)}\right) ds(x)\\
&=0,
\end{aligned}
\en
which leads to the relation
\ben
\left(\int_{\wt D} f(\wt x,k)v_n(\wt x)d\wt x\right)\cdot\left( \int_{H}g(x_d)e^{\sqrt{\lambda_{n}^2-k^2}x_d}dx_d\right)=0 
\enn
for all $n\in \N^{+}$, $k\in[k_{\min},k_{\max}]$.
By Lemma \ref{lem}, there exists a non-empty open subset $K\subset [k_{\min},k_{\max}]$ such that 
\ben
\int_{H}g(x_d)e^{\sqrt{\lambda_n^2-k^2}x_d}dx_d\neq0, \quad k\in K
\enn
holds for all $n\in\mathbb{N}^+$.
For such a non-empty open subset $K\subset [k_{\min},k_{\max}]$, we obtain
\be\label{equ:cn}
\int_{\wt B_R}f(\wt x,k)v_n(\wt x)d \wt x=0, \quad k\in K
\en for all $n\in\mathbb{N}^+$.
This together with the completeness of  
the orthonormal basis $\left\lbrace v_n\right\rbrace _{n=1}^\infty$
 in $L^2(\wt B_R)$ yields
$f(\wt x, k) = 0$ for $\wt x\in \wt B_R$ and $k\in K$. The analyticity of $f$ in $k$ implies the vanishing of $f$ for all $k\in \R$.
\end{proof}

\section{Increasing stability in two dimensions}\label{sec:IS}


In this section we restrict our discussions to the two dimensional settings (i.e., $d=2$) and perform a stability analysis for recovering the wave-number-dependent source $f(x_1,k)$ with explicit dependence on the wavenumber $k\geq1$, which is also well known as the increasing stability with respect to $k$. We retain the notations used in the previous sections. 
For notational convenience, we write $f_k(x_1):=f(x_1,k)$, which is supported in the interval $(-R, R)$.

Introduce the Dirichlet-to-Neumann (DtN) operator $\mathcal{B}: H^{-1/2}(\partial B_R)\rightarrow H^{1/2}(\partial B_R)$ given by $\mathcal{B}u=\frac{\partial u}{\partial\nu}$. Using the DtN operator, we can reformulate the Sommerfeld radiation condition into a transparent
boundary condition
$\frac{\partial u}{\partial\nu}=\mathcal{B}u$ on $\partial B_R$,
where $\nu$ is the unit outer normal on $\partial B_R$. Below we shall derive an upper bound of $\mathcal{B}$ with explicit dependance on the wavenumber, based on the result of \cite{CM2008}. It is worthing mentioning that the estimate for the upper bound also applies to DtN operators  defined on a non-circular closed boundary, although the transparent operator $\mathcal{B}$ is defined on a circle within this paper. 

Let $\widetilde{R}>R$. We denote $D=\{x\in\R^2\setminus \overline{B_R}: |x|<\widetilde{R}\}$ and $V:=\{v\in H^1(D): v=0 \,\, \mbox{on} \,\, |x|=R\}$ equipped with the norm 
\[\|v\|_{V}:=\left(\int_{D}(|\nabla v|^2+k^2|v|^2)dx\right)^{1/2}.\]
 Consider the exterior boundary value problem:
\begin{equation}\label{u2}
\begin{cases}
     \Delta u(x)+k^2u(x)=0, \ \,x\in\R^2\setminus \overline{B_R}, \\
     u=h,  \  \  x\in\partial B_R,\\
    r^{\frac{1}{2}}(\partial_r u-iku)=0, \ \,r:=|x|\rightarrow +\infty,  \\ 
\end{cases}
\end{equation}
where $h\in H^{1/2}(\partial B_R)$. Using the integral equation or variational approach, one can prove that the problem (\ref{u2}) has a unique solution  $u\in H^1_{loc}(\R^2\setminus \overline{B_R})$, even if the boundary $\partial B_R$ is not circular; see e.g., \cite[Chapter 3]{CK2013}. 
By the extension theory in Sobolev spaces,
there exists some $w\in H^1_{loc}(\R^2\setminus B_R)$ satisfying $w|_{\partial B_R}=h$ and $\|w\|_{H^1(D)}\leq C\|h\|_{H^{1/2}(\partial B_R)}$. 
Choose a cut-off function $\chi\in C^\infty(\R^2)$ satisfying $\chi=0$ in $|x|>\frac{R+\widetilde{R}}{2}$ and $\chi=1$ in $|x|\leq R$ .

It is easy to see that the function $\widetilde{u}:=u-\chi w$ satisfies
\begin{equation}\nonumber
\begin{cases}
     \Delta \widetilde{u}+k^2\widetilde{u}=-(\Delta+k^2)(\chi w), \ \,x\in\R^2\setminus \overline{B_R}, \\
     \widetilde{u}=0,  \  \  x\in\partial B_R.\\
\end{cases}
\end{equation}
The right hand side $(\Delta+k^2)(\chi w)$ is compactly supported in $D$,  belonging to the space $H^{-1}(D)$, the dual space of  $H^1_0(D)$.
According to Lemma 3.3 in \cite{CM2008}, we get the following estimate of $\tilde{u}$:
\begin{equation}\label{VR}
\|\widetilde{u}\|_{V}\leq(5+4\sqrt{2}k\widetilde{R})\|(\Delta+k^2)(\chi w)\|_{V^*}.
\end{equation}
Applying the Cauchy-Schwartz inequality and the definition of $V^*$, we deduce that
\begin{equation}\nonumber
\begin{split}
\|\Delta+k^2)(\chi w)\|_{V^*}&=\sup\limits_{\|\varphi\|_{V}\leq1}\langle\Delta+k^2)(\chi w),\varphi  \rangle_{L^2(D)}\\
&=\sup\limits_{\|\varphi\|_{V}\leq1}\int_{D}\big(-\nabla(\chi w)\cdot\nabla \overline{\varphi}+k^2\chi w\cdot\overline{\varphi}\big) dx\\
&\leq Ck^2\|w\|_{H^{1}(D)},
\end{split}
\end{equation}
where $C>0$ depends on $\widetilde{R}$ and $R$ and we have used the fact $k^2\geq1$.
Combining the above inequality and (\ref{VR}), it follows that 
\begin{equation}\nonumber
\begin{split}
\|\widetilde{u}\|_{V}\leq C(5+4\sqrt{2}k\widetilde{R})k^2\|w\|_{H^{1}(D)}\leq C_1k^3\|w\|_{H^{1}(D)}
\end{split}
\end{equation}
for some constant $C_1$ independent of $k$.
Noting that $k\geq1$, we have
\[\|\widetilde{u}\|_{H^{1}(D)}\leq \|\widetilde{u}\|_{V}\leq C_1k^3\|w\|_{H^{1}(D)}.\]
Therefore we obtain
\begin{equation}
\begin{split}
\|u\|_{H^{1}(D)}=\|\widetilde{u}+\chi w\|_{H^{1}(D)}&\leq\|\widetilde{u}\|_{H^{1}(D)}+\|\chi w\|_{H^{1}(D)}\\
&\leq Ck^3\|h\|_{H^{1/2}(\partial B_R)},
\end{split}
\end{equation}
where $C>0$ is independent of $k$.
By the trace theorem, it holds that
\[\|\frac{\partial u}{\partial\nu}\|_{H^{-1/2}(\partial B_R)}\leq C\|u\|_{H^{1}(D)}\]
for some positive constant $C$ depending on $\widetilde{R}$ and $R$. Therefore,
\begin{equation}\label{exp}
\|\mathcal{B}u\|_{H^{-1/2}(\partial B_R)}=\|\frac{\partial u}{\partial\nu}\|_{H^{-1/2}(\partial B_R)}\leq Ck^3\|u\|_{H^{1/2}(\partial B_R)},
\end{equation}
where $C>0$ is independent of $k$. Thus (\ref{exp}) gives an explcit upper bound estimate of $\mathcal{B}$ with respect to $k$.
Numerically, one can also explicitly obtain the Neumann data on $\partial B_R$ once the
Dirichlet date is available on $\partial B_R$, because the the operator $\mathcal{B}$ can be rewritten in a series form.

For the stability estimate we make an additional assumption that there exists a constant $\delta$ such that 
\begin{align}\label{equ:g}
|\widehat{g}(\xi_2)|\geq \delta>0 \quad \mbox{for all} \quad \xi_2\in(-k,k),   
\end{align}
where $\widehat{g}$ denotes the one-dimensional Fourier transform of $g$. Physically, this condition ($\ref{equ:g}$) means that the frequencies of the source function $g$ are mostly restricted to the interval $(-k, k)$.
The main increasing stability result is shown below, which is also valid in three dimensions. 
\begin{thm}\label{1} 
Let $u(x,k)$ be the unique solution to the problem \eqref{equ:uk}-\eqref{equ:SRC}
with $\|f_k\|_{H^1(\R)}\leq M$ for some $M>1$.
Then there exists a constant $C>0$ depending on $R$ and $\delta$ such that 
	\begin{equation}\label{fk}
		\|f_k\|^2_{L^2(\R)}\leq C\big(k^8\epsilon^2+\frac{M^2}{k^{\frac{4}{3}}|\ln \epsilon|^{\frac{1}{2}}} \big), 
	\end{equation}
	where $k\geq1$ and $\epsilon=\|u\|_{H^{1/2}(\partial B_R)}
$.
\end{thm}
We remark that only the Dirichlet boundary data are involved on the right hand of \eqref{fk}.
The proof of Theorem \ref{1} is motivated by the uniqueness proof of Theorem \ref{thm:FT} and the increasing stability for wave-number-independent source functions \cite{CIL2016,LY2017}.

Let $\xi=(\xi_1,\xi_2)\in\R^2$ with $|\xi|=k$.
Multiplying $e^{-i(\xi_1 x_1+\xi_2 x_2)}$ on both sides of (\ref{equ:uk}) and integrating over $ B_R\subset (-R, R)^2$, we obtain (see (\ref{equ:Green1}))
\begin{equation}\label{equ:fg}
	\big(\int_{-R}^Rf(x_1,k)e^{-i\xi_1 x_1}dx_1\big)\big(\int_{-R}^Rg(x_2)e^{-i\xi_2 x_2}dx_2\big)=\int_{\partial B_R}e^{-i\xi\cdot x}(\partial_\nu u+i\xi\cdot\nu u)ds(x),
\end{equation}
Using the estimate of the Nuemann data \eqref{exp}, we can bound the first the firm on the right hand of \eqref{equ:fg} by
\begin{equation}\nonumber
\begin{split}
\int_{\partial B_R}e^{-i\xi\cdot x}(\partial_\nu u) ds(x)&\leq \|\partial_\nu u\|_{H^{-1/2}(\partial B_R)}\|e^{-i\xi\cdot x}\|_{H^{1/2}(\partial B_R)}\\
&\leq\|\partial_\nu u\|_{H^{-1/2}(\partial B_R)}\|e^{-i\xi\cdot x}\|_{H^{1}(\partial B_R)}\\
&\leq Ck^4 \|u\|_{H^{1/2}(\partial B_R)}.
\end{split}
\end{equation}
The second term can be estimated analogously. 
Combining the above inequality and (\ref{equ:g}), (\ref{equ:fg}), we conclude that \begin{equation}
\begin{split}
\Big|\int_{-R}^Rf(x_1,k)e^{-i\xi_1 x_1}dx_1\Big|^2&\leq C\Big|\int_{\partial B_R}e^{-i\xi\cdot x}(\partial_\nu u+i\xi\cdot\nu u)ds(x)\Big|^2\\
&\leq  C\Big|\int_{\partial B_R}e^{-i\xi\cdot x}(\partial_\nu u)ds(x)\Big|^2+C\Big|\int_{\partial B_R}e^{-i\xi\cdot x}(i\xi\cdot\nu u)ds(x)\Big|^2\\
&\leq Ck^8 \|u\|^2_{H^{1/2}(\partial B_R)}.
\end{split}
\end{equation}
Here $C>0$ is independent of $k$.
Hence
\begin{equation}\label{sb}
	\big|\widehat{f_k}(\xi_1)|^2=\big|\frac{1}{\sqrt{2\pi}}\int_{-R}^Rf(x_1,k)e^{-i\xi_1\cdot x_1}dx_1\big|^2\leq Ck^8\epsilon^2 \quad \mbox{for all} \quad |\xi_1|<k.
\end{equation}
Denote
\ben
	I_1(s)=\int_{|\xi_1|\leq s}|\widehat{f_k}(\xi_1)|^2d\xi_1, \quad s>0.
\enn
Let $\xi_1=s\hat{\xi}_1$ for $\hat{\xi}_1\in(-1,1)$. It is easy to get 
\begin{equation}\label{I1s}
		I_1(s)=\int_{|\xi_1|\leq s}|\widehat{f_k}(\xi_1)|^2d\xi_1
		=
		\int_{-1}^1|\widehat{f_k}(s\hat{\xi}_1)|^2s d\hat{\xi}_1\\
		=\frac{1}{2\pi}\int_{-1}^1\big|\int_{-R}^Rf_k(x_1)e^{-is\hat{\xi}_1x_1}dx_1\big|^2s d\hat{\xi}_1.
\end{equation}
Since $e^{-is\hat{\xi}_1x_1}$ is an entire analytic function of $s$, $I_1(s)$  can be extended analytically to  the complex plane with respect to $s$. Thus $I_1(s)=\frac{1}{2\pi}\int_{-1}^1\big|\int_{-R}^Rf_k(x_1)e^{-is\hat{\xi}_1x_1}dx_1\big|^2s d\hat{\xi}_1$ is an entire analytic function of $s=s_1+is_2\in\C$, $s_1,s_2\in\R$ and the following elementary estimates hold.

\begin{lem}\label{I1}
	Let $\|f_k\|_{L^2(\R)}\leq M$.
	Then we have for all $s=s_1+is_2\in\C$ that
	\ben
		|I_1(s)|
		\leq \frac{2RM^2|s|e^{2R|s_2|}}{\pi}.
	\enn
\end{lem}
\begin{proof}
Recall from \eqref{I1s} that
	\begin{equation}
		\nonumber
		I_1(s)=\frac{1}{2\pi}\int_{-1}^1|\int_{-R}^Rf_k(x_1)e^{-is\hat{\xi}_1x_1}dx_1|^2s d\hat{\xi}_1, \quad\mbox{for}\, \ s=s_1+is_2\in\C.
	\end{equation}
	Noting the inequality $|e^{-is\hat{\xi}_1x_1}|\leq e^{R|s_2|}$
	for all $x_1\in(-R,R)$ and $\hat{\xi}_1\in(0,1)$, we obtain
	\begin{equation}
	\begin{split}
		\nonumber
		|I_1(s)|&\leq \frac{1}{2\pi}\int_{-1}^1\big(\int_{-R}^R|f_k(x_1)e^{-is\hat{\xi}_1x_1}|dx_1\big)^2|s| d\hat{\xi}_1\\
		&\leq\frac{1}{2\pi}\int_{-1}^12R\int_{-R}^R|f_k(x_1)|^2|e^{-is\hat{\xi}_1x_1}|^2|s|dx_1 d\hat{\xi}_1
		\leq \frac{2R|s|e^{2R|s_2|}}{\pi}\int_{\R}|f_k|^2dx_1,
	\end{split}
	\end{equation}
	which completes the proof of Lemma \ref{I1} by using the bound of $f_k$.
\end{proof}
Let us recall the following result on the analytical continuation  proved in \cite{CIL2016}.
\begin{lem}\label{2}
	Let $J(z)$ be an analytic function in $S=\{z=x+iy\in \mathbb{C}:-\frac{\pi}{4}<\arg z<\frac{\pi}{4}\}$ and continuous in $\overline{S}$ satisfying
	\begin{equation}
		\nonumber
		\begin{cases}
			|J(z)| \leq \epsilon, \  & z\in (0,L], \\
			|J(z)| \leq V, \  & z\in S,\\
			|J(0)|  =0.
		\end{cases}
	\end{equation}
	Then there exists a function $\mu(z)$ satisfying
	\begin{equation}
		\nonumber
		\begin{cases}
			\mu (z)  \geq\frac{1}{2},\ \ &z\in (L,2^{\frac{1}{4}}L), \\
			\mu (z)  \geq\frac{1}{\pi}((\frac{z}{L})^4-1)^{-\frac{1}{2}},\ \ & z\in (2^{\frac{1}{4}}L, +\infty)
		\end{cases}
	\end{equation}
	such that
	\begin{equation}
		\nonumber
		|J(z)|\leq V\epsilon^{\mu(z)} \quad \mbox{for all} \quad z\in (L, +\infty).
	\end{equation}
\end{lem}
Using Lemma \ref{2}, we show the relation between $I_1(s)$ for $s\in(k,\infty)$ and $I(k)$.
\begin{lem}
	Let $\|f_k\|_{L^2(\R)}\leq M$. Then there exists a function $\mu(s)$ satisfying 
	\begin{equation}
		\nonumber
		\begin{cases}
			\mu (s)  \geq\frac{1}{2},\ \ &s\in (k,2^{\frac{1}{4}}k), \\
			\mu (s)  \geq\frac{1}{\pi}((\frac{s}{k})^4-1)^{-\frac{1}{2}},\ \ & s\in (2^{\frac{1}{4}}k, +\infty)
		\end{cases}
	\end{equation}
	such that 
	\begin{equation}\label{CM}
		|I_1(s)|\leq CM^2k^8e^{(2R+1)s}\epsilon^{2\mu(s)} \quad \mbox{for all}
		\ k<s<+\infty,
	\end{equation}
	where $C>0$ depends on $R$ and $\delta$.
\end{lem}
\begin{proof}
	Let the sector  $S\subset\C$ be given in Lemma \ref{2}. Observe that $|s_2|\leq s_1$ when $s\in S$. It follows from Lemma \ref{I1} that 
	\[|I_1(s)e^{-(2R+1)s}|\leq CM^2,\]
	where $C>0$ depends on $R$.
	Recalling  a prior estimate (\ref{sb}),  we obtain 
	\begin{equation}\label{I-1}
		|I_1(s)|\leq Cs^8\epsilon^2\leq Ck^8\epsilon^2, \quad s\in[0, k],
	\end{equation}
	with  $C>0$ depending on $R$ and $\delta$.
	Hence 
	\ben
		|\frac{1}{k^8}I_1(s)e^{-(2R+1)s}|\leq C\epsilon^2   \quad \mbox{for all} \ s\in[0, k].
	\enn
	Then applying Lemma \ref{2} with $L=k$ to the function $J(s):=\frac{1}{k^8}I_1(s)e^{-(2R+1)s}$,
	we conclude that  there exists a function $\mu(s)$ satisfying 
	\begin{equation}
		\nonumber
		\begin{cases}
			\mu (s)  \geq\frac{1}{2},\ \ &s\in (k, 2^{\frac{1}{4}}k), \\
			\mu (s)  \geq\frac{1}{\pi}((\frac{s}{k})^4-1)^{-\frac{1}{2}},\ \ & s\in (2^{\frac{1}{4}}k, \infty)
		\end{cases}
	\end{equation}
	such that
	\ben
		|\frac{1}{k^8}I_1(s)e^{-(2R+1)s}|\leq CM^2\epsilon^{2\mu(s)},
	\enn
	where $k<s<+\infty$ and $C$ depends on $R$ and $\delta$. Thus we complete the proof.
\end{proof}
Now we show the proof of Theorem \ref{1}.
If $\epsilon\geq e^{-1}$, then we get $M^2e^2\epsilon^2\geq M^2$ and
\[\|f_k\|^2_{L^2(\R)}\leq M^2\leq M^2e^2\epsilon^2\leq M^2e^2k^8\epsilon^2.\]
Hence (\ref{fk}) holds.
If $\epsilon<e^{-1}$, we discuss (\ref{fk}) in two cases as follows.

Case (i):  $2^{\frac{1}{4}}((2R+3)\pi)^{\frac{1}{3}} k^{\frac{1}{3}}< |\ln \epsilon|^{\frac{1}{4}}$. 
Choose $s_0=\frac{1}{((2R+3)\pi)^{\frac{1}{3}}}k^{\frac{2}{3}}|\ln \epsilon|^{\frac{1}{4}}$. It is easy to get $s_0>2^{\frac{1}{4}}k$, which implies
\[-\mu(s_0)\leq-\frac{1}{\pi}((\frac{s_0}{k})^4-1)^{-\frac{1}{2}}\leq-\frac{1}{\pi}(\frac{k}{s_0})^2.\]
A direct application of  estimate (\ref{CM}) shows that 
\begin{eqnarray*}
	|I_1(s_0)|&\leq& CM^2k^8\epsilon^{2\mu(s_0)}e^{(2R+3)s_0}\\
	&\leq& CM^2 k^8 e^{(2R+3)s_0-2\mu(s_0)|\ln\epsilon|}\\
	&\leq& CM^2 k^8 e^{(2R+3)s_0-\frac{2|\ln \epsilon|}{\pi}(\frac{k}{s_0})^2 }\\
	&=& CM^2 k^8 e^{-2(\frac{(2R+3)^2}{\pi})^{\frac{1}{3}}k^{\frac{2}{3}}|\ln \epsilon|^{\frac{1}{2}}(1- \frac{1}{2}|\ln \epsilon|^{-\frac{1}{4}})}.
\end{eqnarray*}
Noting that $1-\frac{1}{2}|\ln \epsilon|^{-\frac{1}{4}}>\frac{1}{2}$ and $(\frac{(2R+3)^2}{\pi})^{\frac{1}{3}}>1$, we have 
\[|I_1(s_0)|\leq CM^2 k^8e^{-k^{\frac{2}{3}}|\ln \epsilon|^{\frac{1}{2}}}.\]
Using the inequality $e^{-t}\leq\frac{14!}{t^{14}}$ for $t>0$, we get
\ben
	|I_1(s_0)|\leq C\frac{M^2 k^8}{(k^{\frac{2}{3}}|\ln \epsilon|^{\frac{1}{2}})^{14}}\leq C\frac{M^2}{k^{\frac{4}{3}}|\ln \epsilon|^{7}}.
\enn
It is clear that $k^{\frac{4}{3}}|\ln \epsilon|^{7} \geq k^{\frac{4}{3}}|\ln \epsilon|^{\frac{1}{2}}$ when $k\geq1$ and $|\ln\epsilon|\geq1$.
Obviously, the following inequalities holds
\ben
	\begin{split}
		\int_{|\xi_1|>s_0}|\widehat{f_k}(\xi_1)|^2d\xi_1&=s_0^{-2}\int_{|\xi_1|>s_0}\xi_1^2|\widehat{f_k}(\xi_1)|^2d\xi_1\\
		&\leq s_0^{-2}\int_{\R}|\widehat{\nabla f_k}(\xi_1)|^2d\xi_1=s_0^{-2}\int_{\R}|\nabla f_k(\xi_1)|^2d\xi_1\leq \frac{M^2}{s_0^2}
	\end{split}
\enn
by the Parseval's identity.
Hence
\begin{equation}\label{ess1}
	\begin{split}
		\|f_k\|^2_{L^2{(\R)}}=\|\widehat{f_k}\|^2_{L^2{(\R)}}&= I_1(s_0)+\int_{|\xi_1|>s_0}|\widehat{f_k}(\xi_1)|^2d\xi_1\\
		&\leq I_1(s_0)+\frac{M^2}{s_0^2}\\
		&\leq C(\frac{ M^2}{k^{\frac{4}{3}}|\ln \epsilon|^{7}}+\frac{M^2}{k^{\frac{4}{3}}|\ln \epsilon|^{\frac{1}{2}}})\\
		&\leq C\frac{ M^2}{k^{\frac{4}{3}}|\ln \epsilon|^{\frac{1}{2}}}.
	\end{split}
\end{equation}

Case (ii):  $|\ln \epsilon|^{\frac{1}{4}}\leq 2^{\frac{1}{4}}((2R+3)\pi)^{\frac{1}{3}} k^{\frac{1}{3}}$. In this case we choose $s_0=k$. Then we find $s_0\geq 2^{-\frac{1}{4}}((2R+3)\pi)^{-\frac{1}{3}}k^{\frac{2}{3}}|\ln \epsilon|^{\frac{1}{4}}$. 
Using  the estimate (\ref{I-1}) gives
$|I(s_0)|\leq Ck\epsilon^2$. Hence
\begin{equation}\label{ess2}
	\begin{split}
		\|f_k\|^2_{L^2{(\R)}}=\|\widehat{f_k}\|^2_{L^2{(\R)}}&= I_1(s_0)+\int_{|\xi_1|>s_0}|\widehat{f_k}(\xi_1)|^2d\xi_1\\
		&\leq I_1(s_0)+\frac{M^2}{s_0^2}\\
		&\leq C(k^8\epsilon^2+\frac{M^2}{k^{\frac{4}{3}}|\ln \epsilon|^{\frac{1}{2}}})\\
		&\leq C(k^8\epsilon^2+\frac{ M^2}{k^{\frac{4}{3}}|\ln \epsilon|^{\frac{1}{2}}}).
	\end{split}
\end{equation}
Combining (\ref{ess1}) and (\ref{ess2}), we finally get
\ben
	\|f_k\|^2_{L^2{(\R)}}\leq C(k^8\epsilon^2+\frac{ M^2}{k^{\frac{4}{3}}|\ln \epsilon|^{\frac{1}{2}}}).
\enn

\section{Numerical examples}\label{sec:numerical}

In this section, we present some numerical results in $\R^2$. Our objective is to demonstrate the feasibility and effectiveness of the reconstruction methods proposed in Section \ref{sec:unique}. The synthetic observation data are generated  by solving the direct problem \eqref{equ:uk}-\eqref{equ:SRC}. 
The radiated multi-frequency data measured on the boundary of a circular region with radius $R$, denoted as
$u(R, \theta; k)$  in polar coordinates, are obtained through the integral expression
\be\label{equ:Ddata}
	 u(R,\theta;k)=\int_{\B_R}\Phi (x-y;k) f(y_1,k)g(y_{2}) dy, \quad \theta \in (0,2\pi], \, k\in [k_{\min},k_{\max}].
\en
The Neumann data is subsequently computed using the Dirichlet-to-Neumann map,  given by

\be\label{equ:Ndata}
        \partial_{\nu}u(R,\theta; k)=\sum_{n=-\infty}^{+\infty}k\frac{H^{(1)\prime}_{n}(kR)}{H^{(1)}_{n}(kR)}u_{n}(R,k)e^{in\theta},	\quad \theta \in (0,2\pi], \, k\in [k_{\min},k_{\max}],
\en
where $H_n^{(1)}(\cdot)$ denotes the first kind Hankel function of order $n$, $H_n^{(1)\prime}(\cdot)$ is the first derivative of $H_n^{(1)}(\cdot)$ and $u_{n}(R, k)$ is calculated using the Fourier series representation
\ben
u_{n}(R,k)=\frac{1}{2\pi}\int_{0}^{2\pi}u(R,\theta; k)e^{-in\theta}d\theta, \quad k\in [k_{\min},k_{\max}].
\enn
The noise polluted data $u^{\delta} (R, \theta; k)$ are generated by
 \be \label{equ:noiseData}
 u^{\delta}	(R,\theta; k)=u(R,\theta; k)+\delta \zeta |u(R,\theta; k)|,\quad \theta \in (0,2\pi],\, k\in [k_{\min},k_{\max}],
 \en 
where $\zeta$ is a uniformly distributed random number within the range $[-1, 1]$ and $\delta$ represents the noise level. The Neumann data with noise can be computed similarly by
\ben\label{noiseNdata}
\partial_{\nu}u^{\delta}(R,\theta; k)=\sum^{+\infty}_{n=-\infty}k\dfrac{H^{(1)\prime}_{n}(kR)}{H^{(1)}_{n}(kR)}u^\delta_{n}(R,k)e^{in\theta},\quad \theta \in (0,2\pi], k\in [k_{\min},\, k_{\max}],
\enn
where $u^\delta_{n}(R, k)$ denote the Fourier coefficients of 
$u^\delta(R,\theta; k)$. 

%

\begin{figure}[H]
	\centering 
	\includegraphics[scale=0.3]{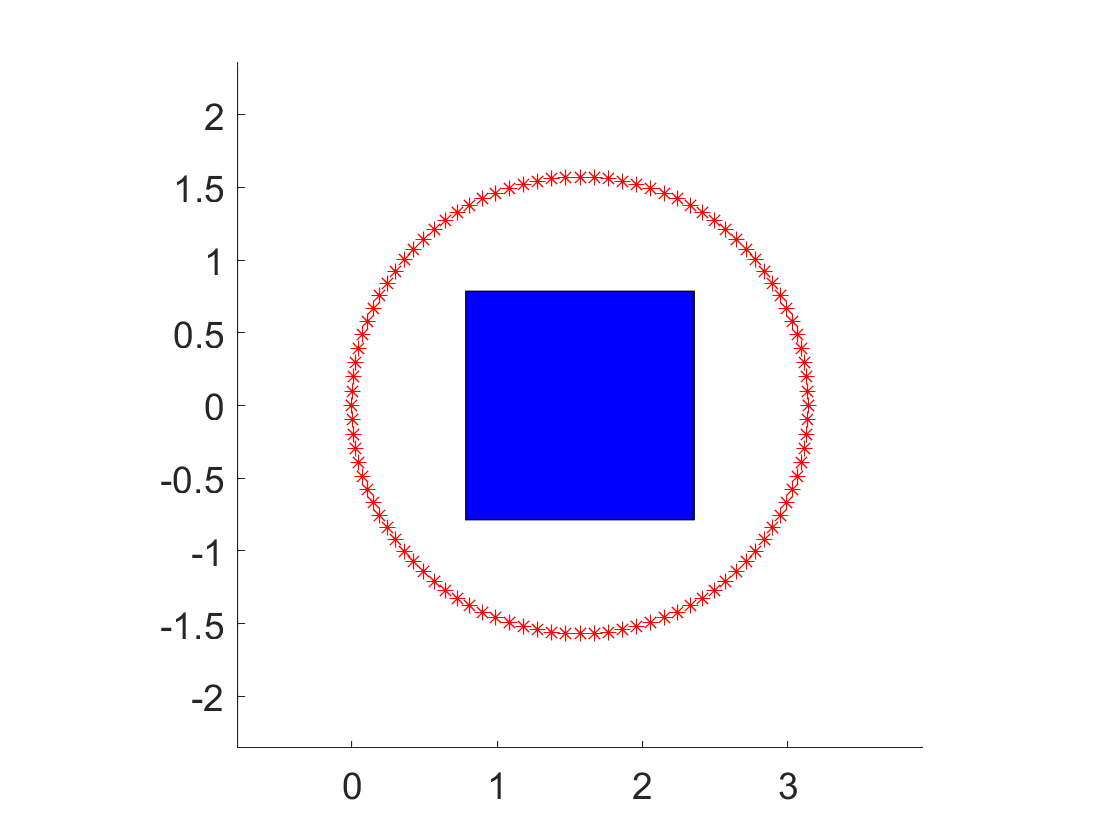} 
	\caption{The blue rectangle represents the source support, where $\supp f(\cdot,k)=[\pi/4,3\pi/4]$ and $\supp g(\cdot)=[-\pi/4,\pi/4]$. The red arsterisks show the measurement points on the circle centered at $(\pi/2,0)$ with radius of $\pi/2$.} 
	\label{fig:PS}
\end{figure}

Now we provide details of our numerical implementation and carry out a series of experiments in two dimensions. In our work, we set the lower bound of the wave number as $k_{\min}=0$ and the upper bound as $k_{\max}=K$. The radiated field data are collected as follows:
\be\label{data:D}
\lbrace u(R,\theta_m; k_j):\, \theta_m=2\pi m/M,\, k_j=jK/J,\, m=1,2,\ldots, M, \, j=1,2,\ldots, J\rbrace.
\en
The data are measured on the circle $\partial \B_R(x_0)$, where $\partial \B_R(x_0)$ represents the circle centered at $x_0=(\pi/2,0)$ with the radius of $R=\pi/2$. Additionally, we calculate the Neumann data by using the formula defined in \eqref{equ:Ndata} as follows:
\be\label{data:N}
\lbrace \partial_\nu u(R,\theta_m; k_j):\, \theta_m=2\pi m/M,\, k_j=jK/J,\, m=1,2,\ldots, M, \, j=1,2,\ldots, J\rbrace.
\en
Unless otherwise specified, we always set $M=100$ in our experiments. 
In Fig.\ref{fig:PS}, we show the support of the source $f(x_1,k)g(x_2)$ for $k\in[k_{\min},k_{\max}]$ along with the depiction of the measurement boundary $\partial \B_R(x_0)$.
For the sake of simplicity, we assume that $g(x_{2})=1$ if $x_2\in [-\pi/4,\pi/4]$ and $g=0$ if else.

\subsection{Numerical implementation by the Dirichlet-Laplacian method}


In this subsection, by using the expansion defined in \eqref{equ:expansion}, we reconstruct the exact source $f(\cdot,k) \in L^2(0,\pi)$ for some fixed $k \in [0,K]$.  To achieve this, one can opt for the eigenfunctions of the Dirichlet Laplacian, as they form a complete and orthonormal set in the space $L^2(0,\pi)$. Consider the one-dimensional eigenvalue problem 
\be\label{equ:DLFpi}
	\left\{
	\begin{array}{l}
		\begin{aligned}
			-v^{''}(x_{1})&=\lambda^2 v(x_{1}), \quad x_{1}\in[0,\pi],\\
			v(x_{1}) &= 0, \qquad\qquad x_{1}=0, \, \pi.
		\end{aligned}
	\end{array}
	\right. 
\en
We obtain the set of Dirichlet eigenvalues as $\lbrace n^2:\,|n|\in\mathbb{N}^+\rbrace$ and the set of Dirichlet eigenfunctions as $\lbrace \sin(nx_{1}):\,|n|\in\mathbb{N}^+\rbrace$. 
To perform the expansion given in \eqref{equ:expansion}, we choose a suitably large value of $N$ such that $|n|\leq N$. By substituting the eigenvalues and eigenfunctions calculated from \eqref{equ:DLFpi} into \eqref{equ:Green2}, i,e., $\varphi_n=\sin(nx_{1})e^{\sqrt{n^{2}-k^{2}}x_{2}}$ for $|n|\leq N$ in \eqref{equ:Green2}, one can deduce that
\ben
	\begin{aligned}
		&\int_{\B_{R}}f(x_{1},k)g(x_{2})\sin(nx_{1})e^{\sqrt{n^{2}-k^{2}}x_{2}}dx\\
		&=\int_{\partial \B_R}\left\lbrace \partial_{\nu}u(x,k)\sin(nx_{1})e^{\sqrt{n^{2}-k^{2}}x_{2}}-u(x,k)\partial_{\nu}\left(\sin(nx_{1})e^{\sqrt{n^{2}-k^{2}}x_{2}}\right) \right\rbrace ds(x).
	\end{aligned}	
\enn
Note that $B_R=B_R(x_0)$ with $x_0=(\pi/2, 0)$ in our experiments.
It implies that
\be\label{equ:fN}	
   \int_{0}^{\pi}f(x_1,k_j)\sin(nx_{1})dx_{1}=\dfrac{M_{n,j}}{G_{n,j}}	
\en
for $n=\pm1,\pm2,\cdots,\pm N$ and $j=1,2,\cdots,J$, where 
\ben
  M_{n,j}:=\int_{\partial \B_R}\left\lbrace \partial_{\nu}u(x,k_j)\sin(nx_{1})e^{\sqrt{n^{2}-k_j^{2}}x_{2}}-u(x,k_j)\partial_{\nu}\left(\sin(nx_{1})e^{\sqrt{n^{2}-k_j^{2}}x_{2}}\right) \right\rbrace ds(x),
\enn
\ben
  G_{n,j}:={\int_{-\frac{\pi}{4}}^{\frac{\pi}{4}}g(x_{2})e^{\sqrt{n^{2}-k_j^{2}}x_{2}}dx_{2}}.	
\enn
We approximate $M_{n,j}$ through the discrete measurement data in polar coordinates as defined in \eqref{data:D} and \eqref{data:N}. 
The reconstructed source $f(x_1,k_j)$ can be approximated by  $f_N(x_{1},k_j)$ with $x_1\in [0,\pi]$ and $k_j\in[0,K]$.
Here $f_N(x_{1},k_j)$
can be expressed as 
\be\label{equ:expansionN}
	f_N(x_{1},k_j)=2\sum_{n=1}^N\,f_{n,j}\sin(nx_{1}),\quad  f_{n,j}=-f_{-n,j}.
\en 
According to \eqref{equ:fN} and \eqref{equ:expansionN}, we deduce coefficients $f_{n,j}$ from \be\label{equ:coefficientN}
f_{n,j}=\dfrac{1}{\pi}\int_{0}^{\pi}f(x_{1},k_j)\sin(nx_{1})dx_1=\dfrac{1}{\pi}\dfrac{M_{n,j}}{G_{n,j}},\quad n=1,2,\cdots,N, \, j=1,2,\cdots,J.
\en
It is worthy noting that here we have assumed $G_{n,j}\neq 0$ for all $n=1,2,\cdots,N$ and $j=1,2,\cdots, J$.
Similarly, the reconstruction formula using noise data is given by 
\be\label{equ:noisereconstructN}
f^{\delta}_N(x_{1},k_j)=2\sum_{n=1}^N\,f^{\delta}_{n,j}\sin(nx_{1}),\quad x_1\in [0,\pi],\, k_j\in[0,K],\\
f^{\delta}_{n,j}=\dfrac{1}{\pi}\dfrac{M^{\delta}_{n,j}}{G_{n,j}},\quad n=1,2,\cdots,N, \, j=1,2,\cdots,J,\label{equ:noisecoefficients}
\en 
where 
\ben
M^{\delta}_{n,j}:=\int_{\partial \B_R}\left\lbrace \partial_{\nu}u^{\delta}(x,k_j)\sin(nx_{1})e^{\sqrt{n^{2}-k_j^{2}}x_{2}}-u^{\delta}(x,k_j)\partial_{\nu}\left(\sin(nx_{1})e^{\sqrt{n^{2}-k_j^{2}}x_{2}}\right) \right\rbrace ds(x).
\enn
The above reconstruction method is referred to the Dirichlet-Laplacian method.

\begin{figure}[H]
	\centering	
	\subfigure[$f_1(x_1,k)$, $k=0.5$, $N=17$]{
		\includegraphics[scale=0.25]{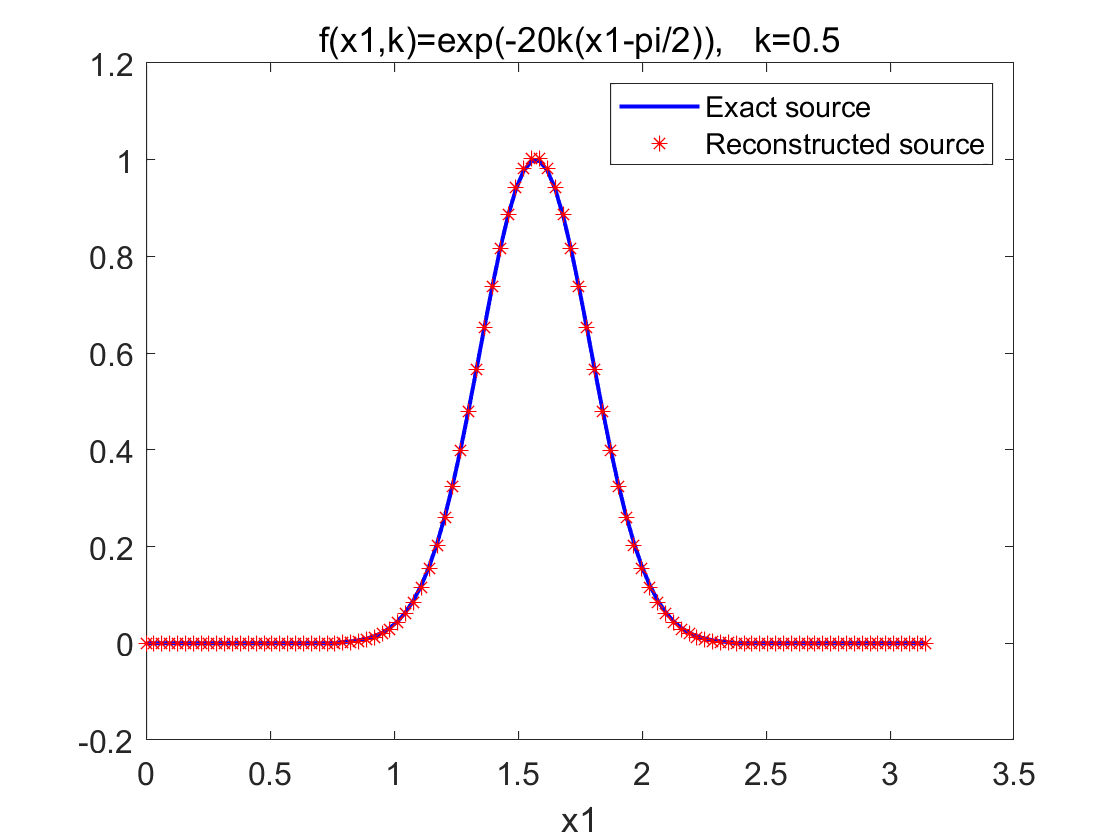}}
	\subfigure[$f_2(x_1,k)$, $k=0.5$, $N=25$]{
		\includegraphics[scale=0.25]{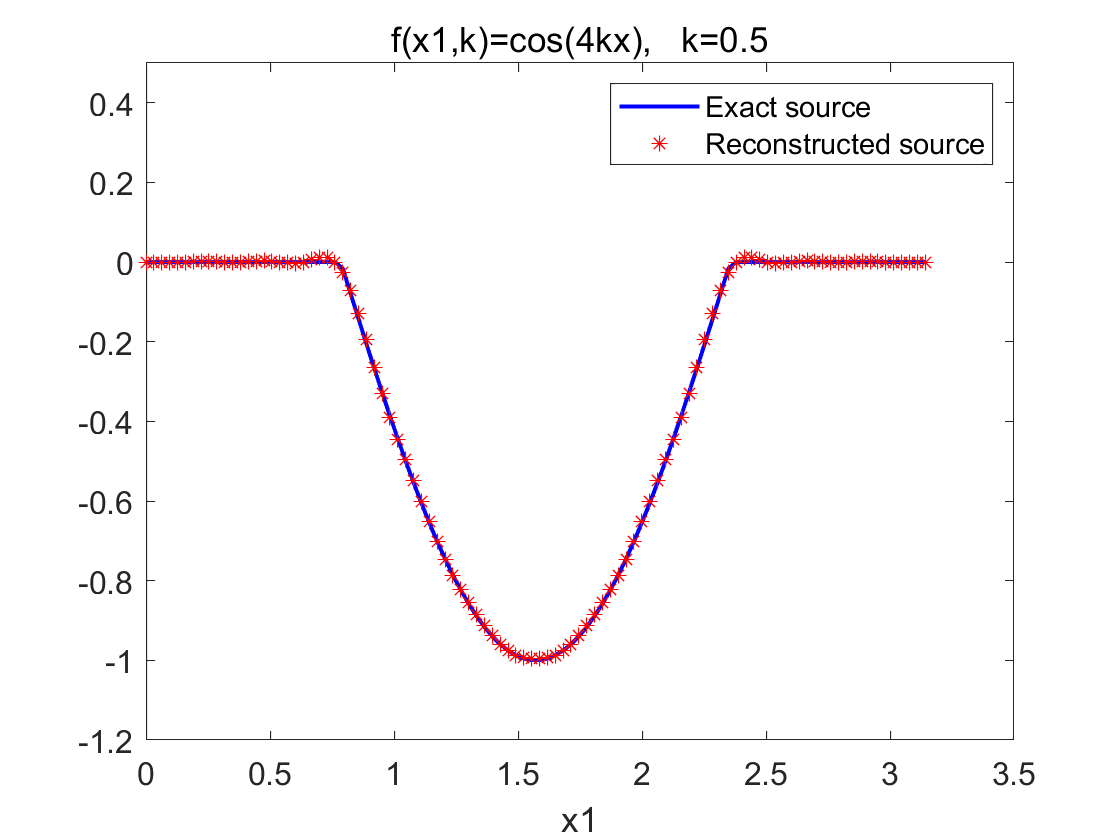}}
	\subfigure[$f_3(x_1,k)$, $k=0.5$, $N=26$]{
		\includegraphics[scale=0.25]{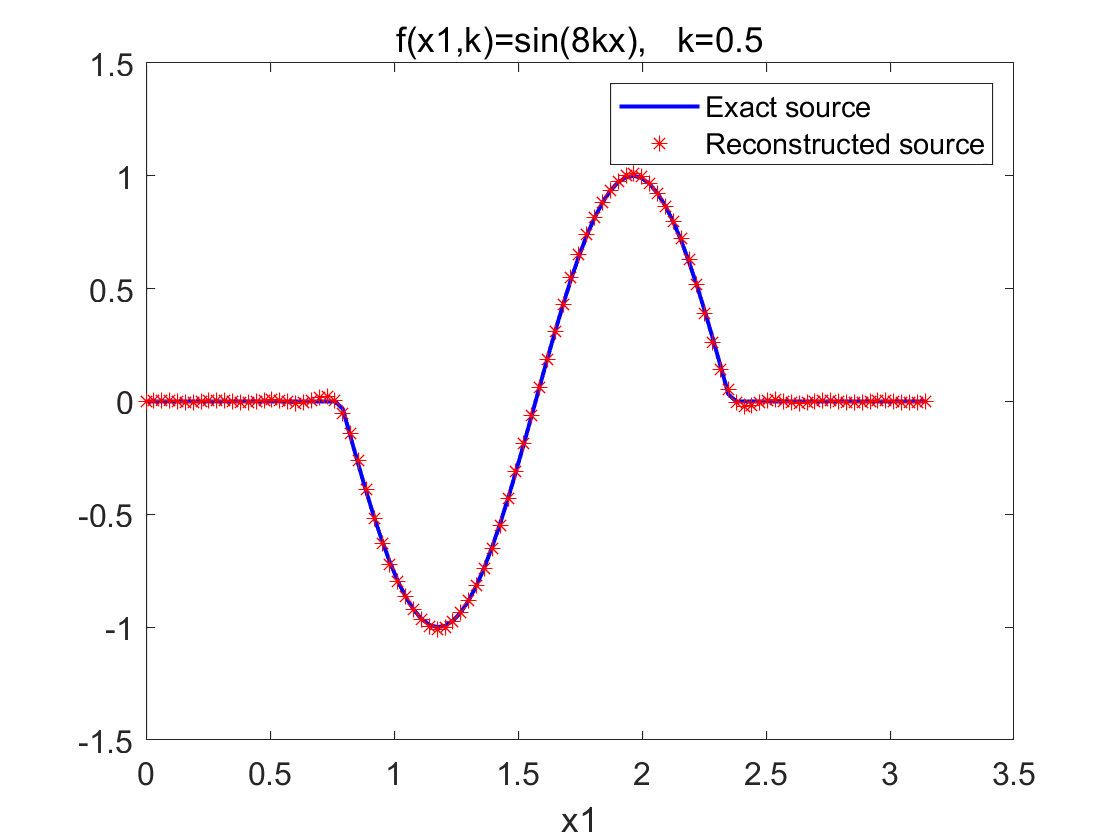}}
	\subfigure[$f_4(x_1,k)$, $k=0.5$, $N=26$]{
		\includegraphics[scale=0.25]{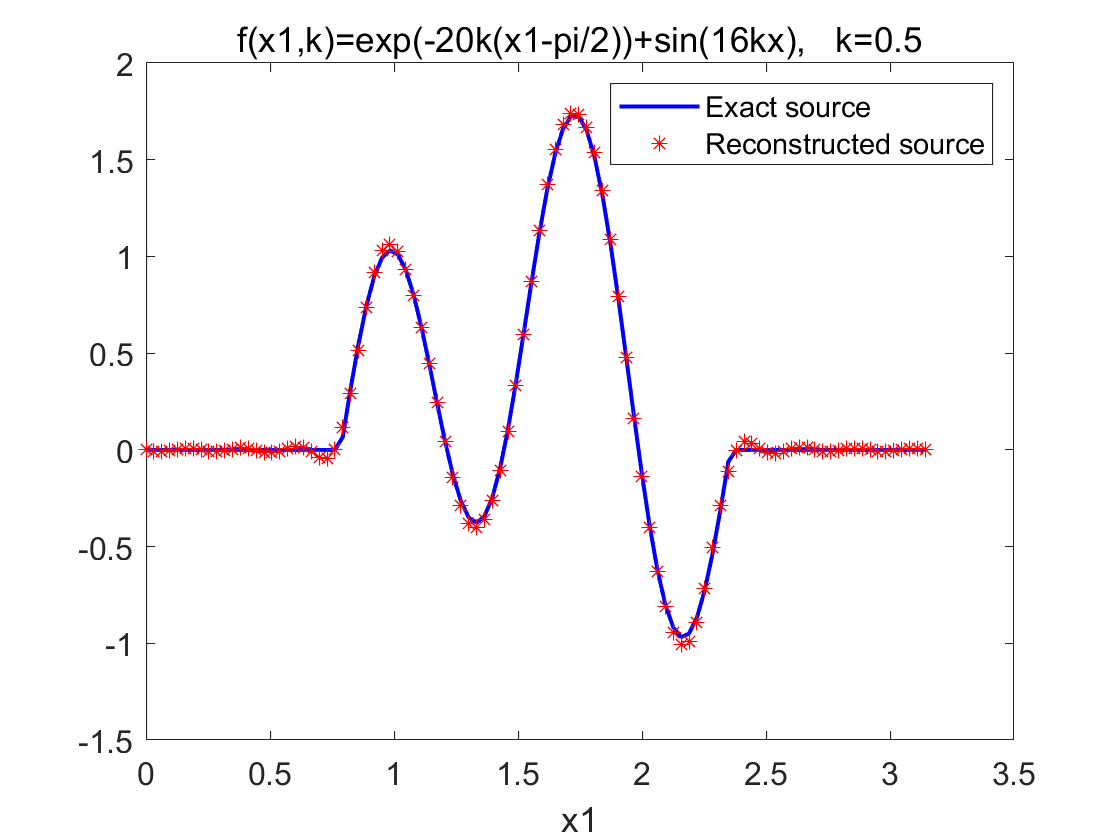}}		
	\caption{The reconstructed source by the Dirichlet- Laplacian method with a fix $k=0.5$. The four source functions $f_j(x_1,k)$ ($j=1,2,3,4$) shown in (a)-(d) are reconstructed with $N=17,\, 25,\, 26,\,26$, respectively. }
	\label{fig:DL:k05}
\end{figure}

\par Using formulas \eqref{equ:expansionN} and \eqref{equ:coefficientN}, we demonstrate the reconstructions of the following source functions when $k=0.5$: 
\be\label{equ:f1}
f_1(x_{1},k)=\left\{
\begin{array}{l}
	e^{-20k(x_{1}-\pi/2)^{2}},\quad x_{1}\in[\frac{\pi}{4},\frac{3\pi}{4}],\\
	0,\qquad\qquad \qquad\quad   {\rm otherwise},		
\end{array}
\right.
\en
\be\label{equ:f2}
f_2(x_{1},k)=\left\{
\begin{array}{l}
	\cos(4kx_{1}),\quad  x_{1}\in[\frac{\pi}{4},\frac{3\pi}{4}],\\
	0,\qquad\qquad\quad  {\rm otherwise},		
\end{array}
\right.
\en
\ben\label{equ:f3}
f_3(x_{1},k)=\left\{
\begin{array}{l}
	\sin(8kx_{1}),\quad x_{1}\in[\frac{\pi}{4},\frac{3\pi}{4}],\\
	0,\qquad\qquad\quad   {\rm otherwise},		
\end{array}
\right.
\enn
\ben\label{equ:f4}
f_4(x_{1},k)=\left\{
\begin{array}{l}
	e^{-20k(x_{1}-\pi/2)^{2}}+\sin(8kx_{1}),\quad x_{1}\in[\frac{\pi}{4},\frac{3\pi}{4}],\\
	0,\qquad\qquad \qquad\qquad \qquad\qquad {\rm otherwise}.		
\end{array}
\right.
\enn
It is obvious that the source can be well reconstructed by choosing an appropriate parameter $N$ in Fig.\ref{fig:DL:k05}. We use the source function $f_1(x_1,0.5)$ to generate the noise-free radiation field data on $\partial \B_R$. The reconstructed source $f_N(x_1,0.5)$ is displayed in Fig.\ref{fig:DL:k05} (a) with $N=17$. In Fig.\ref{fig:DL:k05} (b), the source function $f_2(x_1,0.5)$ can be well reconstructed with $N=25$. Choosing the source function $f_3$ or $f_4$ with $k=0.5$, we obtain the reconstruction results in Fig.\ref{fig:DL:k05} (c) and Fig.\ref{fig:DL:k05} (d) with $N=26$, respectively. 

 
However, the above method fails in the case $G_{n,j}=0$, because the denominator on the right side of \eqref{equ:coefficientN} vanishes. For example, in reconstructing the source function $f_1(x_1,k)$ with $k=1$, or the source function
\ben
f_5(x_{1},k)=\left\{
\begin{array}{l}
	\cos(2/3kx_1),\quad x_{1}\in[\frac{\pi}{4},\frac{3\pi}{4}],\\
	0,\qquad\qquad \qquad {\rm otherwise},		
\end{array}
\right.
\enn
with $k=3$, one cannot obtain a correct reconstruction with $N=10$. 

\begin{figure}[H]
	\centering	
	\subfigure[$k_1=0.99$]{
		\includegraphics[scale=0.25]{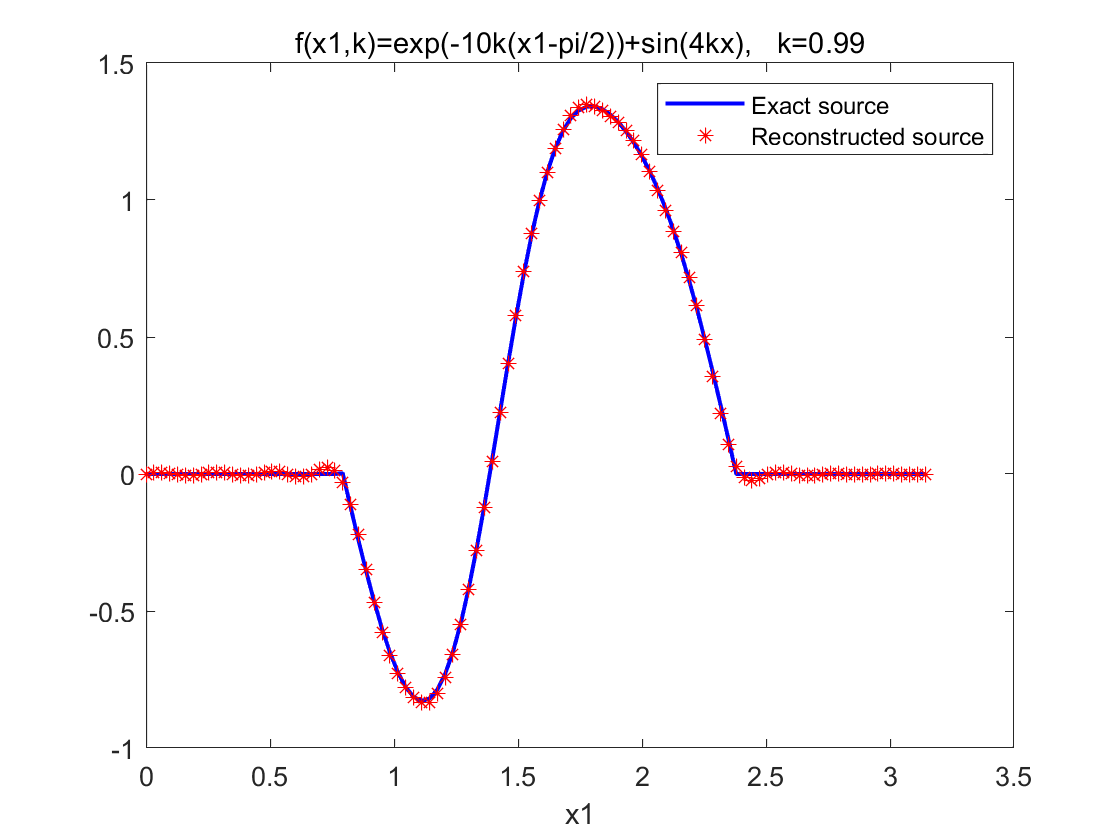}}
	\subfigure[$k_2=1.99$]{
		\includegraphics[scale=0.25]{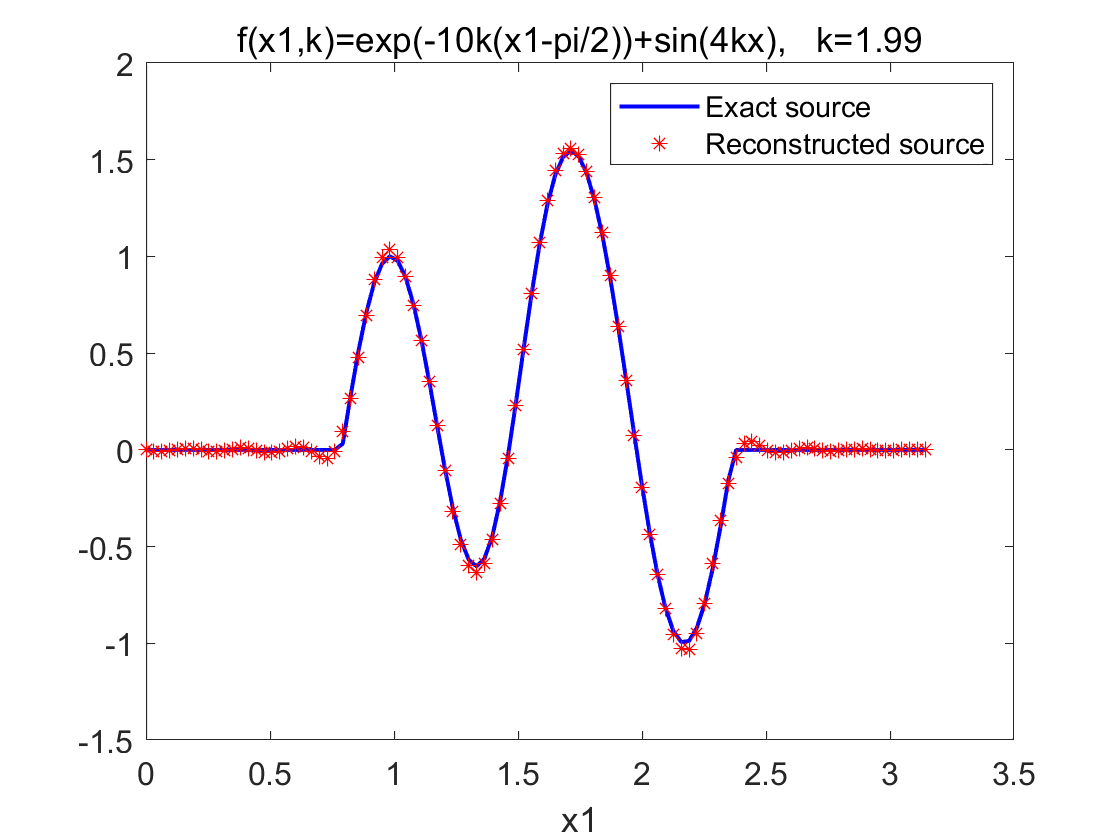}}
	\subfigure[$k_3=2.99$]{
		\includegraphics[scale=0.25]{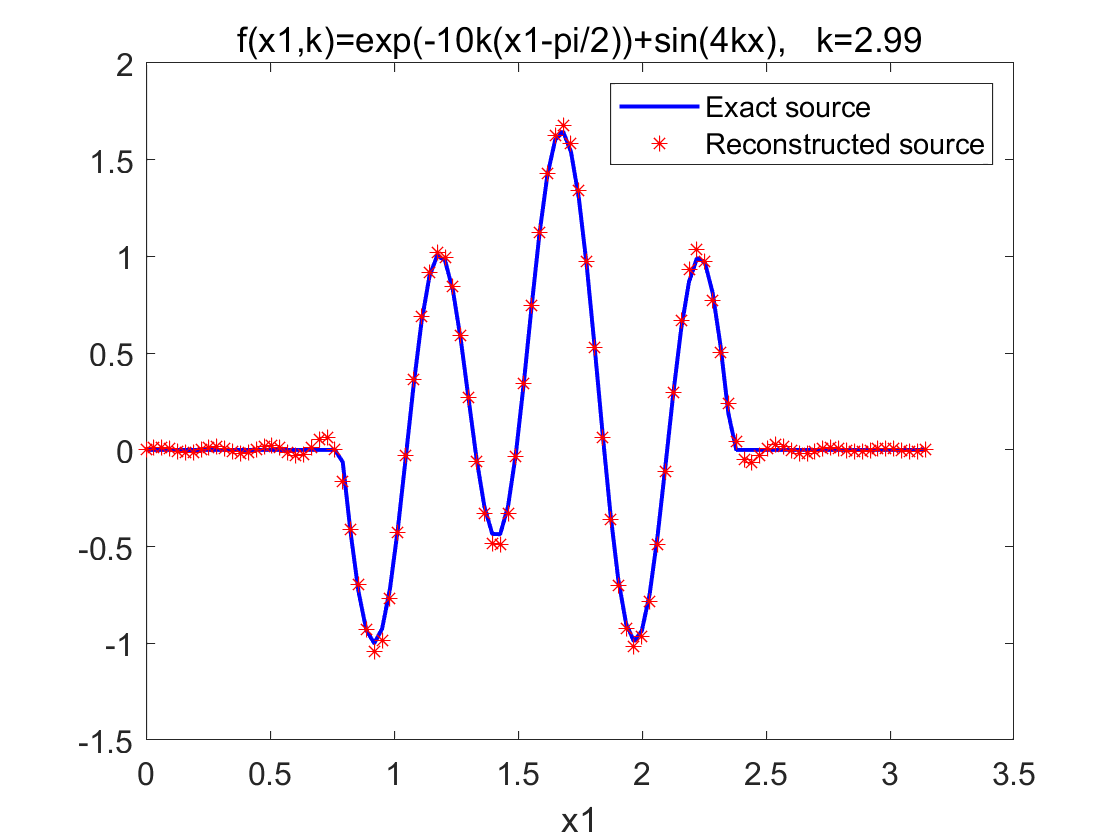}}
	\subfigure[$k_4=3.99$]{
		\includegraphics[scale=0.25]{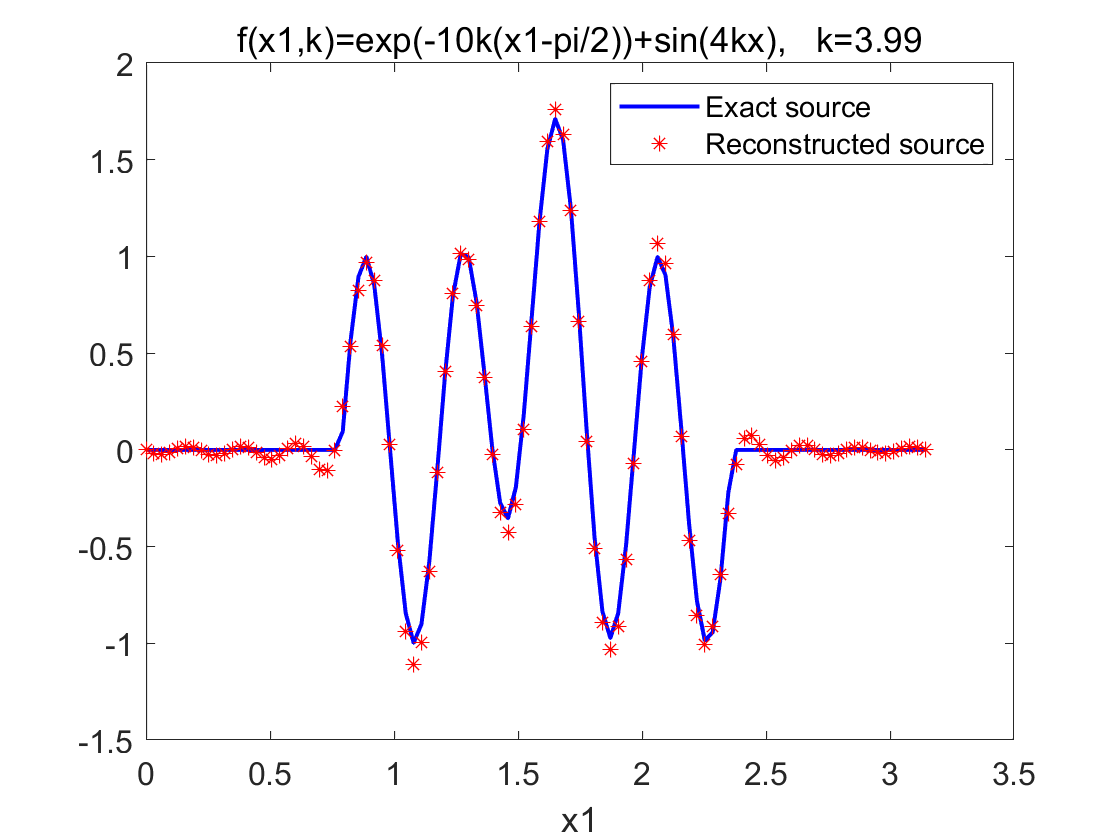}}		
	\caption{The reconstructed results for the source function $x_1\mapsto e^{-10k_j(x_1-\pi/2)^2}+\sin(4k_jx_1)$ supported in $[\pi/4,3\pi/4]$ and $N=26$ by the Dirichlet-Laplacian method. Here we use four different wave-numbers $k_1=0.99$, $k_2=1.99$, $k_3=2.99$ and $k_4=3.99$.}
	\label{fig:DL:S4}
\end{figure}

Next, we choose different wave-numbers $0<k\notin\mathbb{N}^+$ to test the effectiveness of the Dirichlet-Laplacian method. Fig.\ref{fig:DL:S4} illustrates that the source function $e^{-10k(x_{1}-\pi/2)^{2}}+\sin(4kx_{1})$ supported in $[\pi/4,3\pi/4]$ can be well reconstructed for various wave-numbers $0<k\notin\mathbb{N}^+$. When $k_1=0.99$ and $N=26$ (as shown in Fig. \ref{fig:DL:S4} (a)), it is evident that the reconstructed source closely matches the exact source. In other cases, such as $k_2=1.99$, $k_3=2.99$ and $k_4=3.99$, we also achieve satisfactory reconstructions with $N=26$.

\begin{figure}[H]
	\centering
	\subfigure[ $\delta=0.5\%$, $N=6$, ${\rm error}=0.1015$]{
		\includegraphics[scale=0.25]{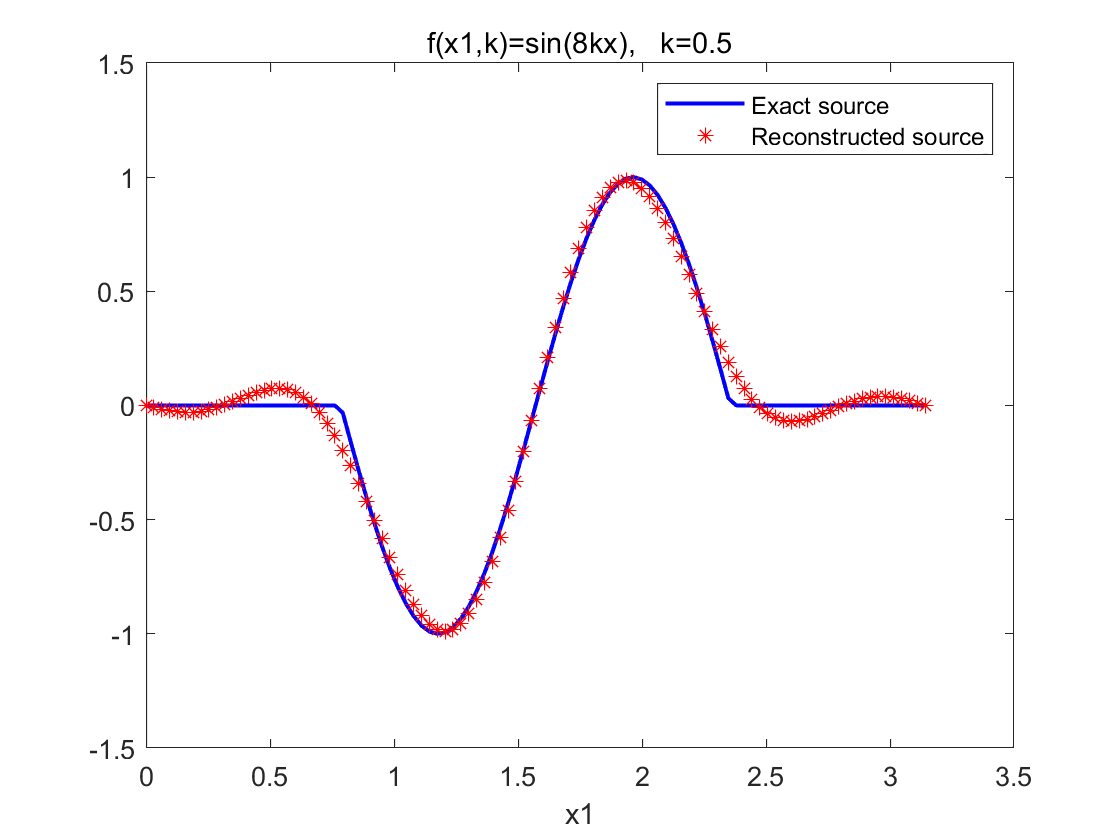}}	
	\subfigure[ $\delta=2\%$, $N=6$, ${\rm error}=0.2226$]{
		\includegraphics[scale=0.25]{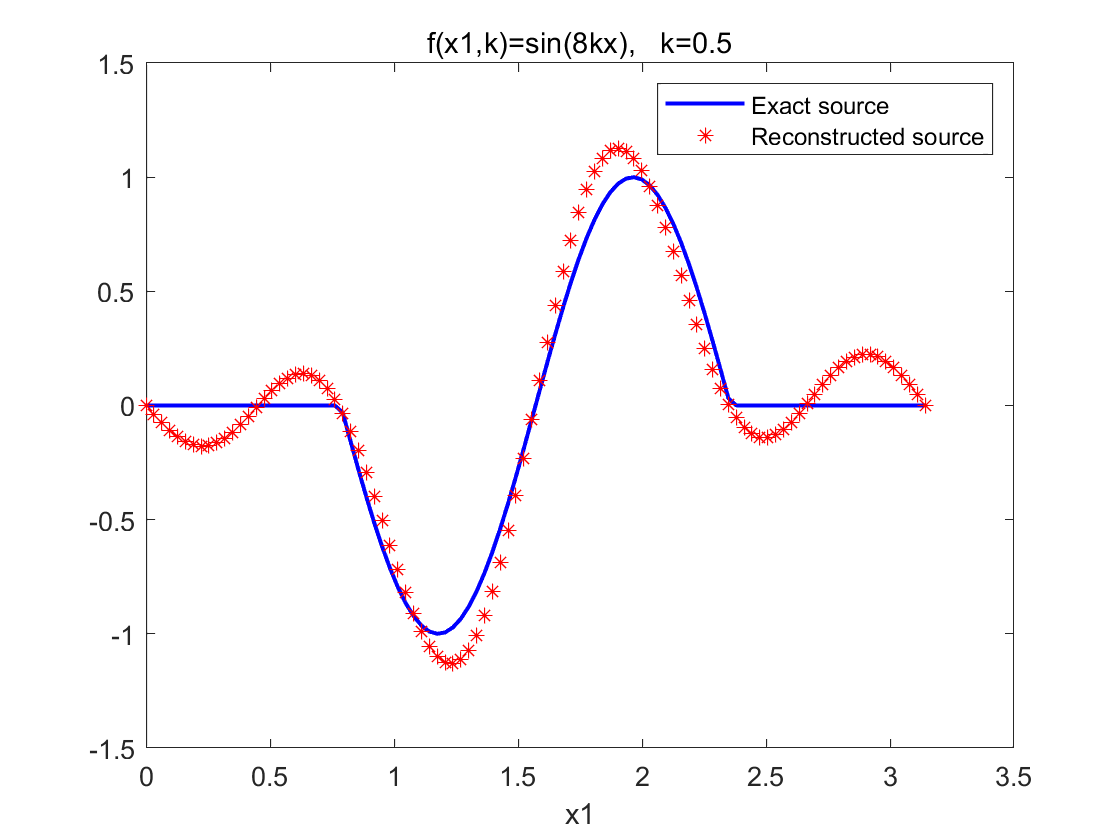}}
	\subfigure[ $\delta=10\%$, $N=4$, ${\rm error}=0.3739$]{
		\includegraphics[scale=0.25]{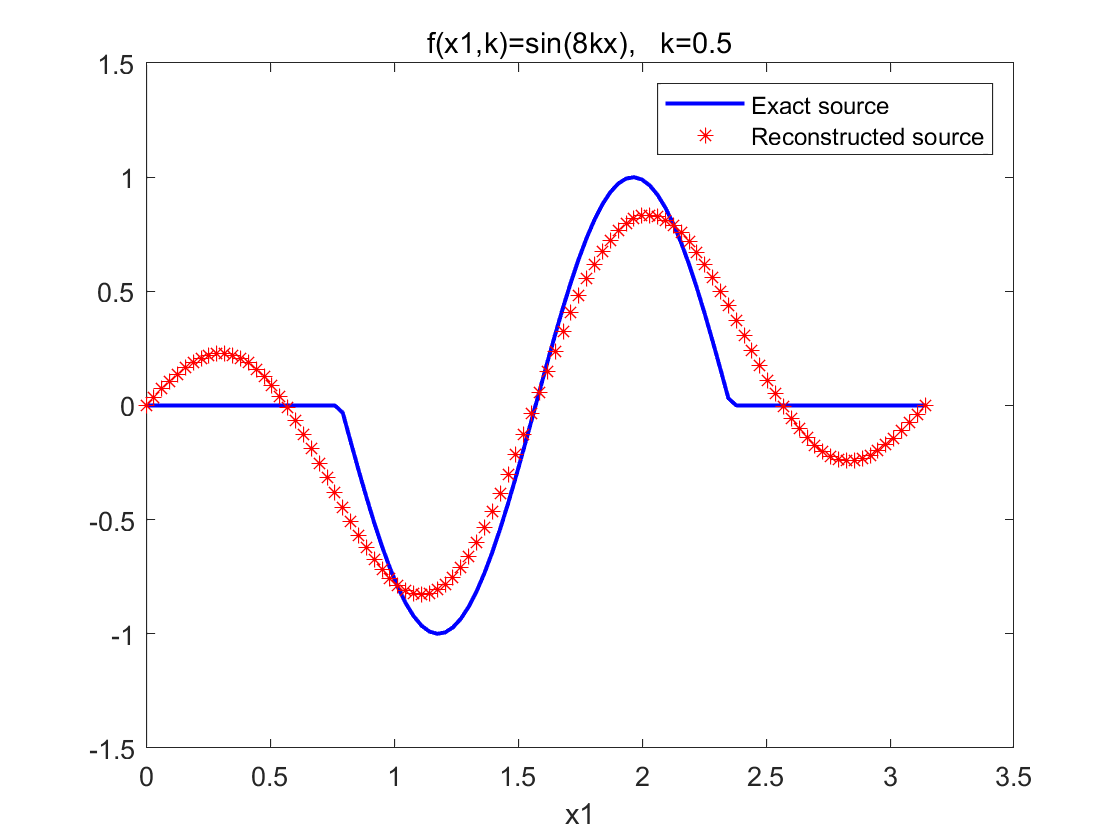}}
	\subfigure[ $\delta=30\%$, $N=4$, ${\rm error}=0.4383$]{
		\includegraphics[scale=0.25]{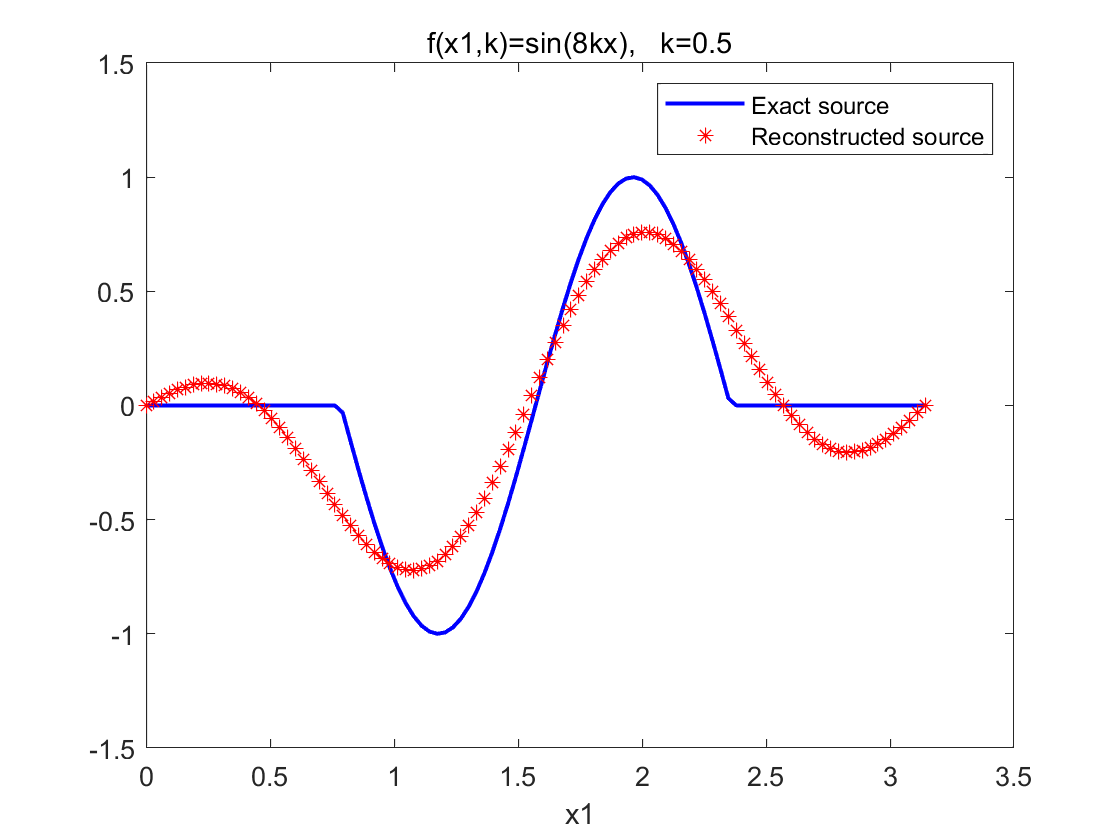}}		
	\caption{The reconstructed source by the Dirichlet-Laplacian method with different noise levels $\delta$. }
	\label{fig:DL:noise}
\end{figure}

Finally, we examine the sensitivity of the Dirichlet-Laplacain method to random noise. The measurement data \eqref{equ:Ddata} is polluted by random noise using the formula \eqref{equ:noiseData}. 
Fig.\ref{fig:DL:noise} displays the reconstructions $f_N^\delta$ of the source function $f=\sin(8kx_1)$ supported in $[\pi/4,3\pi/4]$ with noise levels $\delta=0.5\%$, $\delta=2\%$, $\delta=10\%$ and $\delta=30\%$. It can be concluded that the error decreases as the noise level decreases, where the error is defined as
\be\label{equ:error}
{\rm error}=\dfrac{\| f-f_N^\delta \|_2}{\|f\|_2}.
\en

%


\subsection{Numerical Implementation by the Fourier-Transform Method}
The purpose of this subsection is to reconstruct the exact source $f(\cdot,k)\in L^2(0,\pi)$ for some $k\in[0,K]$ by using the Fourier expansion. A computational formula for each Fourier coefficient is also established. We still cosider the source function considered as in Fig.\ref{fig:PS}. According to \eqref{equ:test1}, the Fourier basis functions in $L^2[0,\pi]$ are given by 
\ben
e^{i\xi_1x_1}, \quad \xi_1^2+\xi_2^2=k^2,\,k\in[0,K],
\enn 
where $\xi_1=2n$, $n\in \mathbb{Z}$, since  $\xi_1=\frac{2\pi}{2R}n=2n$ by our choice $R=\pi/2$. Next, we establish computational formulas for the Fourier coefficients of $f(x_1,k)\in L^2[0,\pi]$ under the Fourier basis functions $\{e^{i(2n)x_1}: x_1\in [0,\pi],\, n\in\mathbb{Z}\}$. According to formula (\ref{equ:Green1}), we obtain
\ben
\begin{aligned}
	&\int_{\B_{R}}f(x_{1},k)g(x_{2})e^{-i(\xi_1x_1+\xi_2x_2)}dx\\
	&=\int_{\partial \B_R}\{ \partial_{\nu}u(x,k)+i\xi \cdot \nu u(x,k)\} e^{-i(\xi_1x_1+\xi_2x_2)} ds(x).
\end{aligned}	
\enn
A simple calculation yields that
\be\label{equ:tidlefN}	
\int_{0}^{\pi}f(x_{1},k_j)e^{-i (2n) x_1}dx_{1}=\dfrac{\tilde{M}_{n,j}}{\tilde{G}_{n,j}}
\en
for $n\in \mathbb{Z}$, $j=1,2,\cdots, W$, where $\tilde{M}_{n,j}$ and $\tilde{G}_{n,j}$ are defined as
\ben
\tilde{M}_{n,j}:=
\int_{\partial \B_R}\{ \partial_{\nu}u(x,k_j)+i(2n,\sqrt{k_j^2-4n^2}) \cdot \nu u(x,k_j)\} e^{-i(2n x_1+\sqrt{k_j^2-4n^2} x_2)}  ds(x),
\enn
\ben
\tilde{G}_{n,j}:={\int_{-\frac{\pi}{4}}^{\frac{\pi}{4}}g(x_{2})e^{-i\sqrt{k_j^2-4n^2}x_{2}}dx_{2}}.
\enn
The reconstructed source $\tilde{f}_N(x_{1},k_j)$ with $x_1\in [0,\pi]$, $k_j\in[0,K]$ can be written as 
\be\label{equ:expansiontidleN}
\tilde{f}_N(x_{1},k_j)=\sum_{\left| n\right| \leq N}\,\tilde{f}_{n,j}e^{i (2n) x_1}.
\en
From \eqref{equ:tidlefN}, the coefficients are approximately computed by
\be\label{equ:coefficienttidleN}
\tilde{f}_{n,j}\approx\dfrac{1}{\pi}\int_{0}^{\pi}f(x_{1},k_j)e^{-i (2n) x_1}dx_{1}=\dfrac{1}{\pi}\dfrac{\tilde{M}_{n,j}}{\tilde{G}_{n,j}},\quad n\in \mathbb{Z}, \, j=1,2,\cdots,J,
\en
where the denominator is again assumed to be non-vanishing.
Similarly, the reconstruction formula from noise data is given by 
\ben
\tilde{f}^{\delta}_N(x_{1},k_j)=\sum_{\left| n\right| \leq N}\,\tilde{f}_{n,j}e^{i (2n) x_1},\quad x_1\in [0,\pi],\, k_j\in[0,K],\\
\tilde{f}^{\delta}_{n,j}\approx\dfrac{1}{\pi}\dfrac{\tilde{M}^{\delta}_{n,j}}{\tilde{G}_{n,j}},\quad n\in \mathbb{Z}, \, j=1,2,\cdots,W,
\enn 
where 
\ben
\tilde{M}^{\delta}_{n,j}:=\int_{\partial \B_R}\{ \partial_{\nu}u^\delta(x,k_j)+i(2n,\sqrt{k_j^2-4n^2}) \cdot \nu u^\delta(x,k_j)\} e^{-i(2n x_1+\sqrt{k_j^2-4n^2} x_2)}  ds(x).
\enn
These reconstruction formulas are quite similar to the Dirichlet-Laplacian method.

\par
\begin{figure}[H]
	\centering	
	\subfigure[$e^{-20k(x_1-\pi/2)^2}$, $N=9$]{
		\includegraphics[scale=0.25]{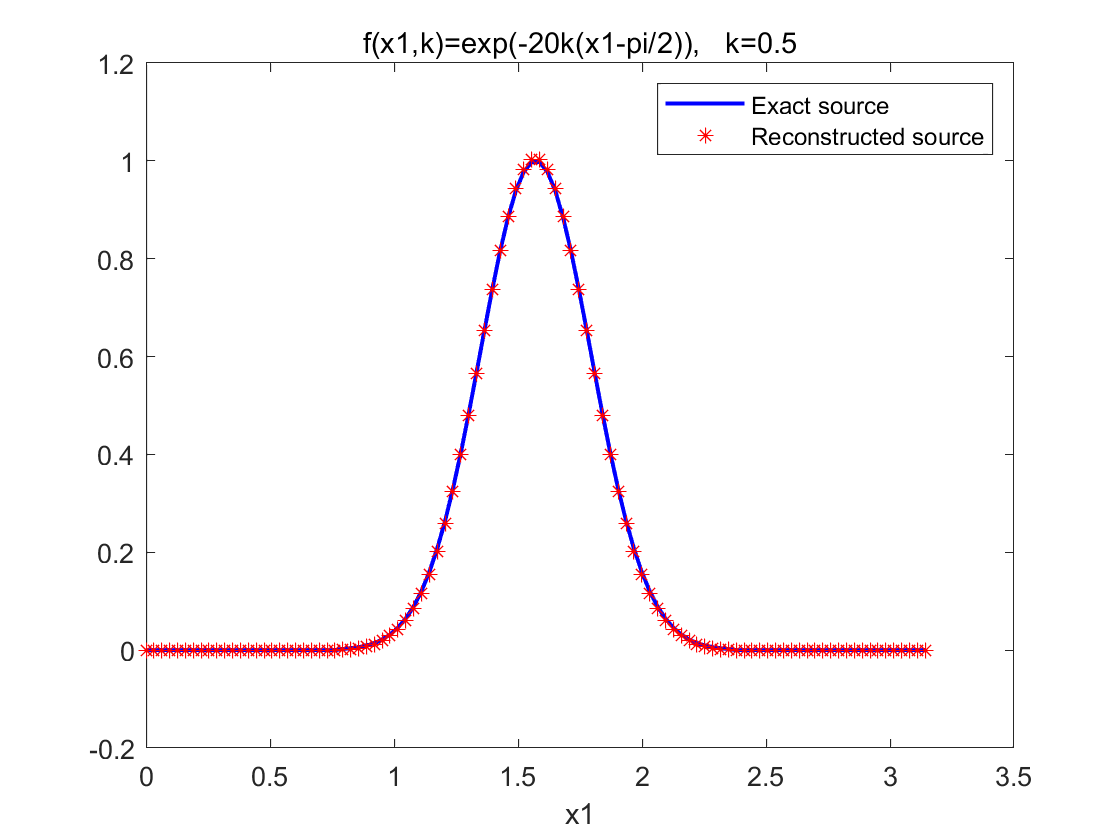}}
	\subfigure[$\cos(4kx_1)$, $N=12$]{
		\includegraphics[scale=0.25]{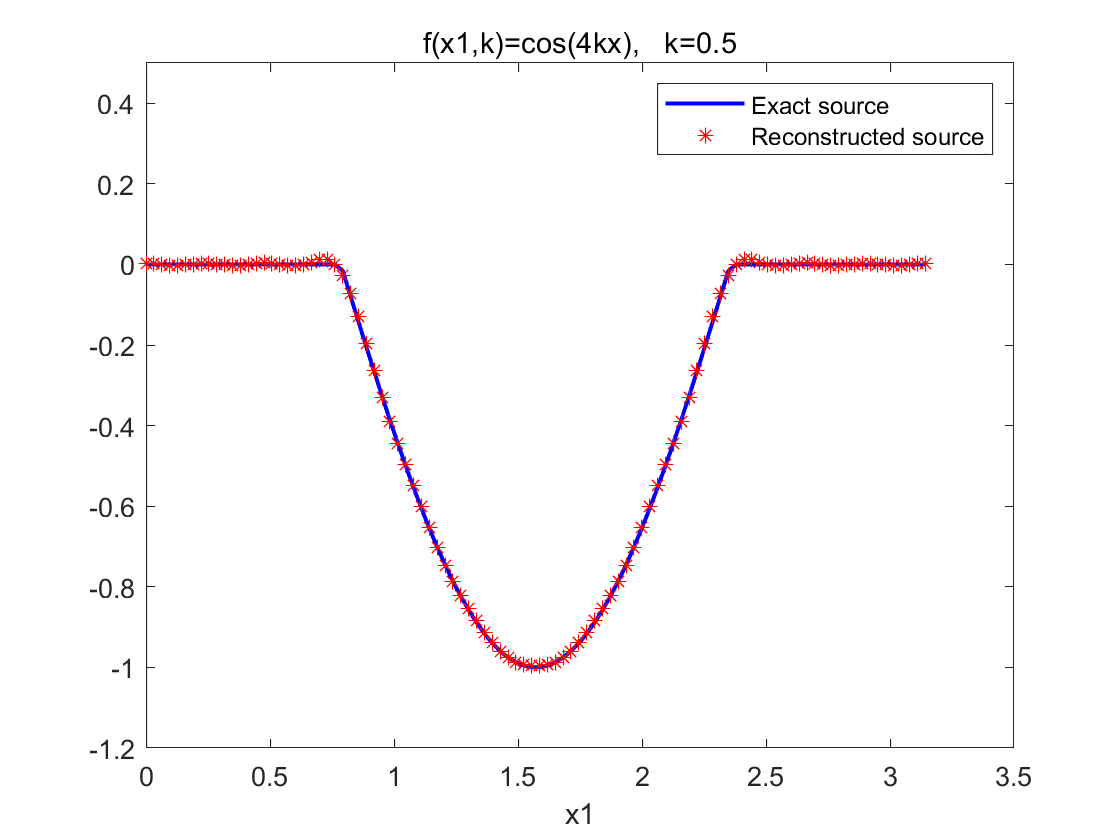}}
	\subfigure[$\cos(12kx_1)$, $N=12$]{
		\includegraphics[scale=0.25]{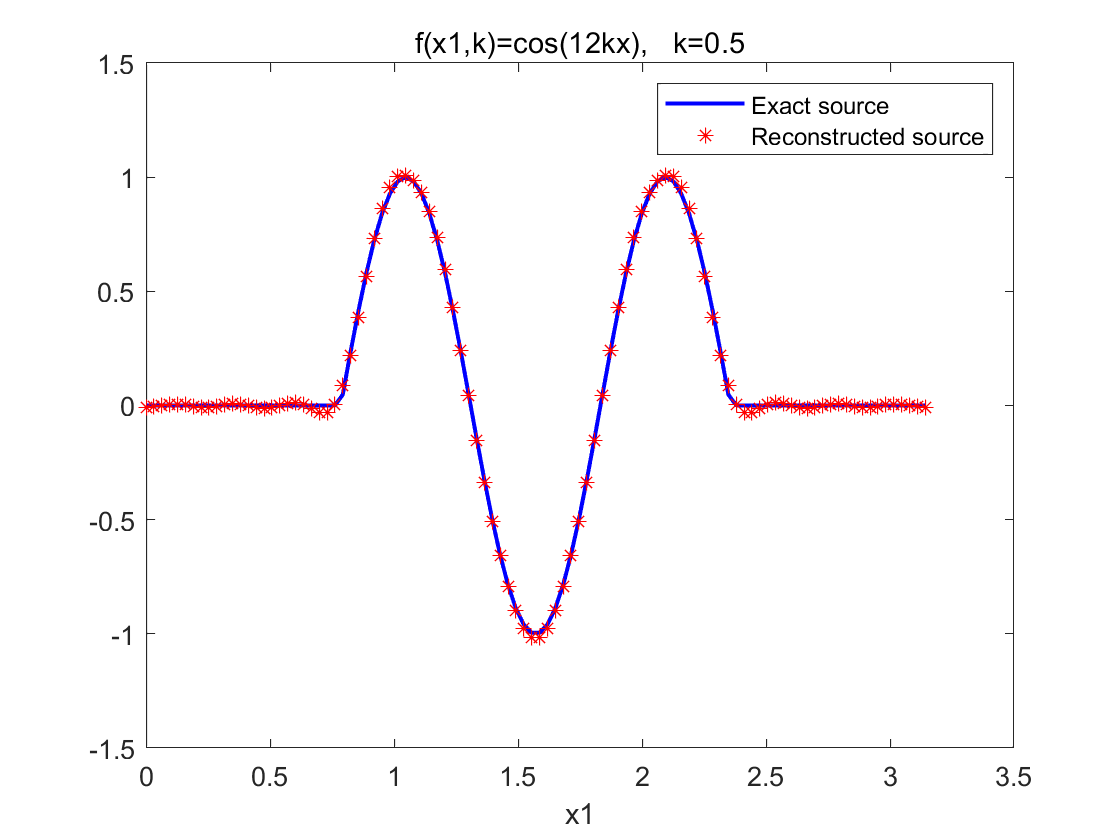}}
	\subfigure[$e^{-20k(x_1-\pi/2)^2}+\cos(12kx_1)$, $N=12$]{
		\includegraphics[scale=0.25]{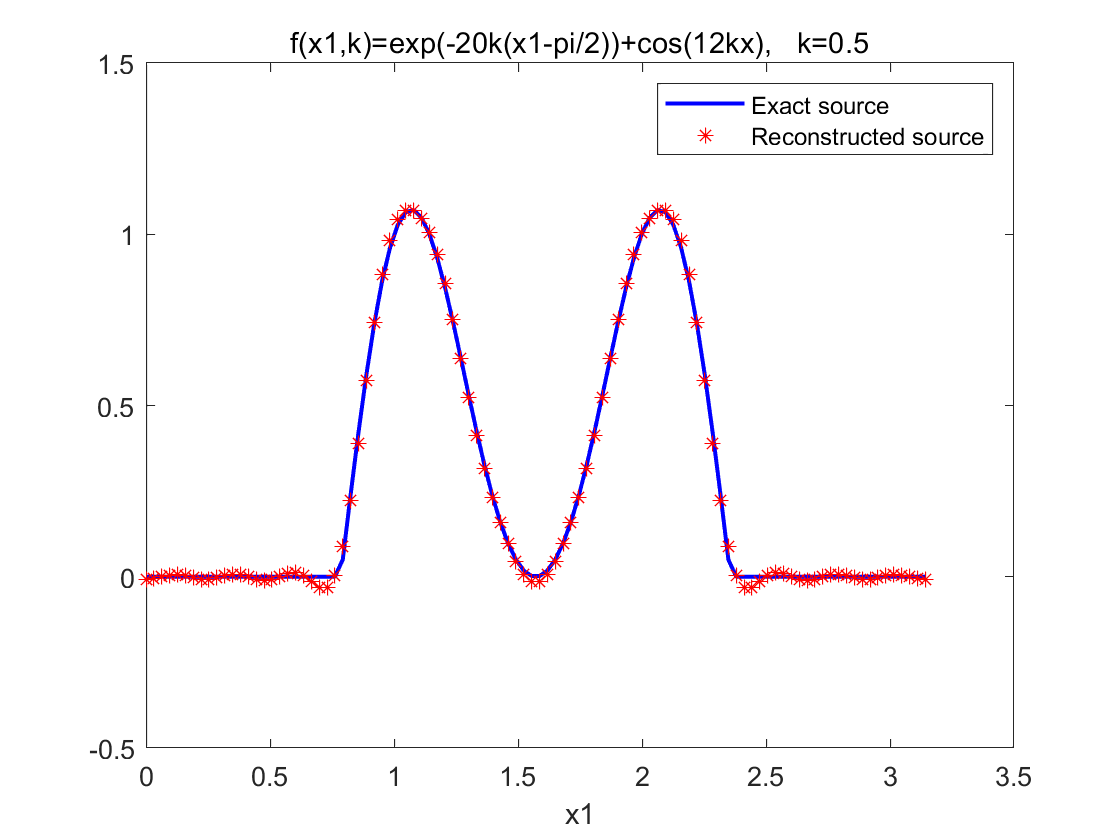}}		
	\caption{The reconstructed source by the Fourier-Transform method with $k=0.5$. Four different sources are reconstructed with $N=9,\,12,\,12,\,12$, respectively.
	}
	\label{fig:FT:k05}
\end{figure}

 Using formulas \eqref{equ:expansiontidleN} and \eqref{equ:coefficienttidleN}, we present the reconstructions  when $k=0.5$. Choosing the source function $f_1(x_1,k)$ in \eqref{equ:f1}, the reconstructed results are shown in Fig.\ref{fig:FT:k05} (a) with $N=9$. In Fig.\ref{fig:FT:k05} (b), the source function $f_2(x_1,k)$, identical to $\eqref{equ:f2}$, is well reconstructed with $N=12$ and $k=0.5$. For source functions $\cos(12kx_{1})$ supported in $[\pi/4,3\pi/4]$ with $k=0.5$ and $e^{-20k(x_{1}-\pi/2)^{2}}+\cos(12kx_{1})$ supported in $[\pi/4,3\pi/4]$ with $k=0.5$, the experimental results shown in Fig.\ref{fig:FT:k05} (c) with $N=12$ and Fig.\ref{fig:DL:k05} (d) with $N=12$ demonstrate that the reconstructed sources closely resemble the exact one.

\begin{figure}[H]
	\centering	
	\subfigure[$\cos(kx_1)$, $k=2$, $N=12$]{
		\includegraphics[scale=0.25]{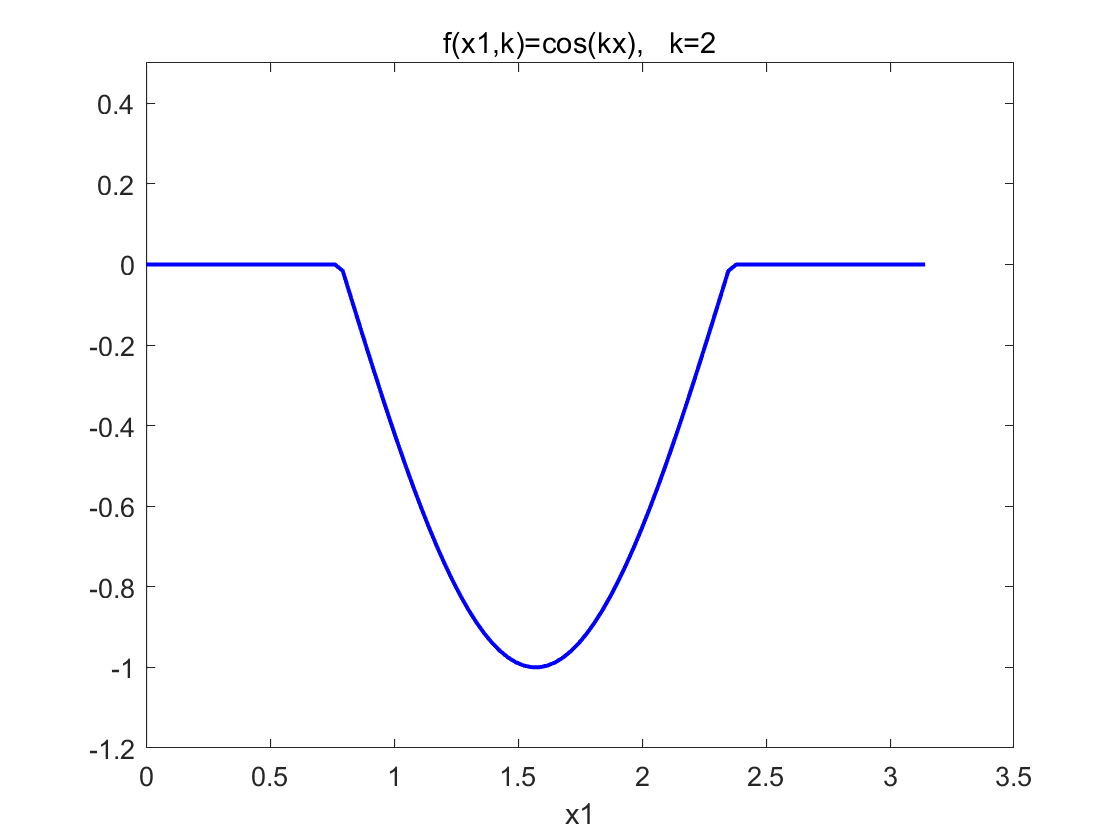}}
	\subfigure[$\cos(kx_1)$, $k=2.001$, $N=12$]{
		\includegraphics[scale=0.25]{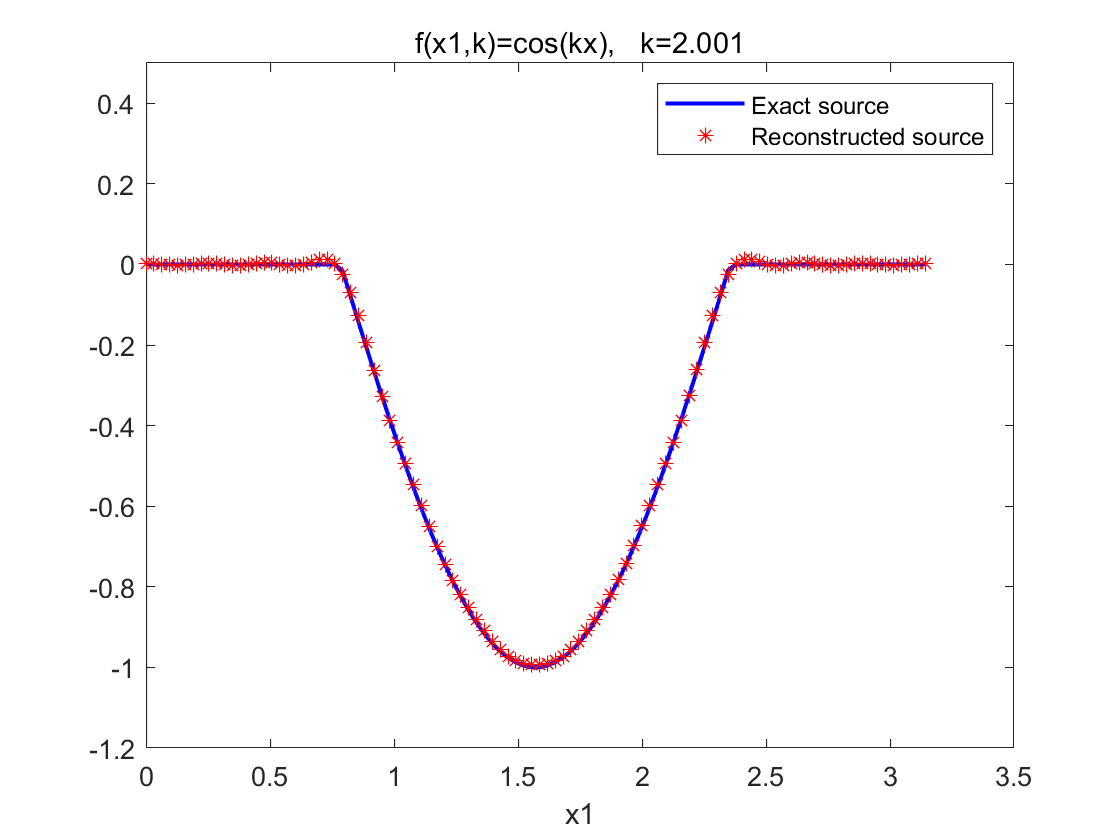}}	
	\caption{The reconstructed source by the Fourier-Transform method with $N=12$. The left shows the exact source function with $k=2$, while the right shows the reconstruction result with $k=2.001$.}
	\label{fig:FT:k}
\end{figure}

\begin{figure}[H]
	\centering	
	\subfigure[$k=1.5$, $N=8$]{
		\includegraphics[scale=0.23]{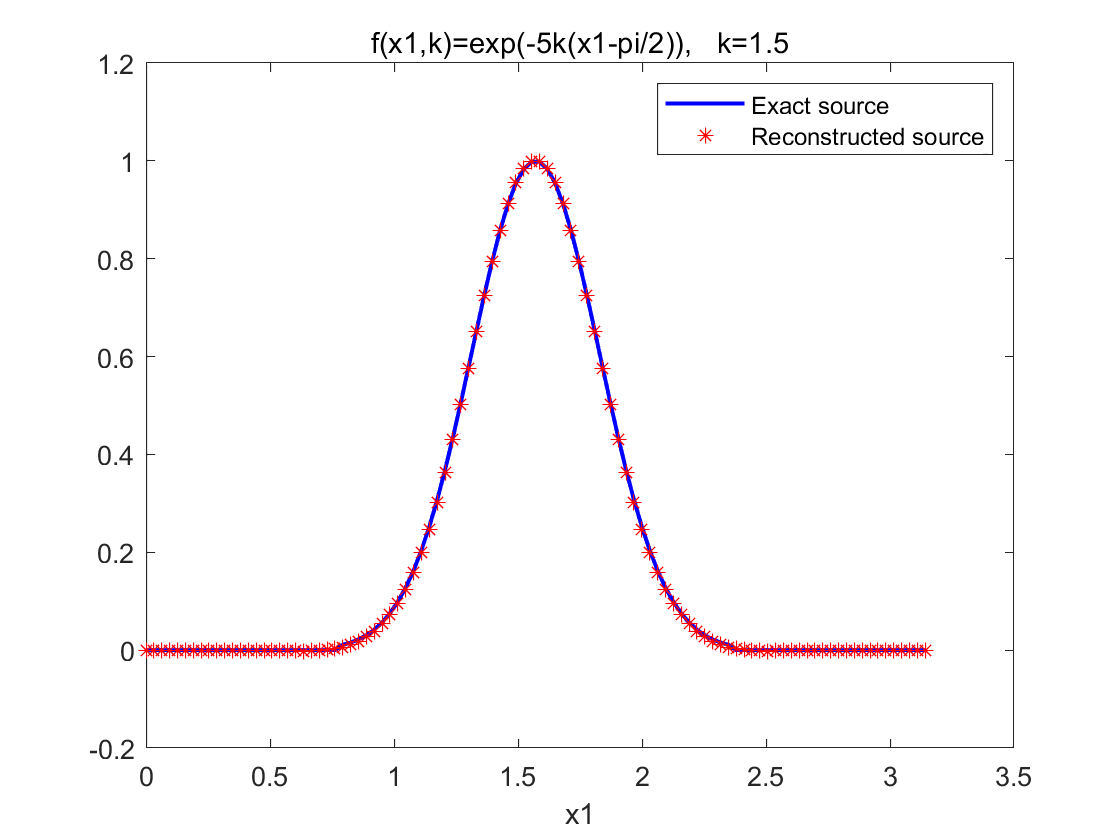}}
	\subfigure[$k=3.5$, $N=8$]{
		\includegraphics[scale=0.23]{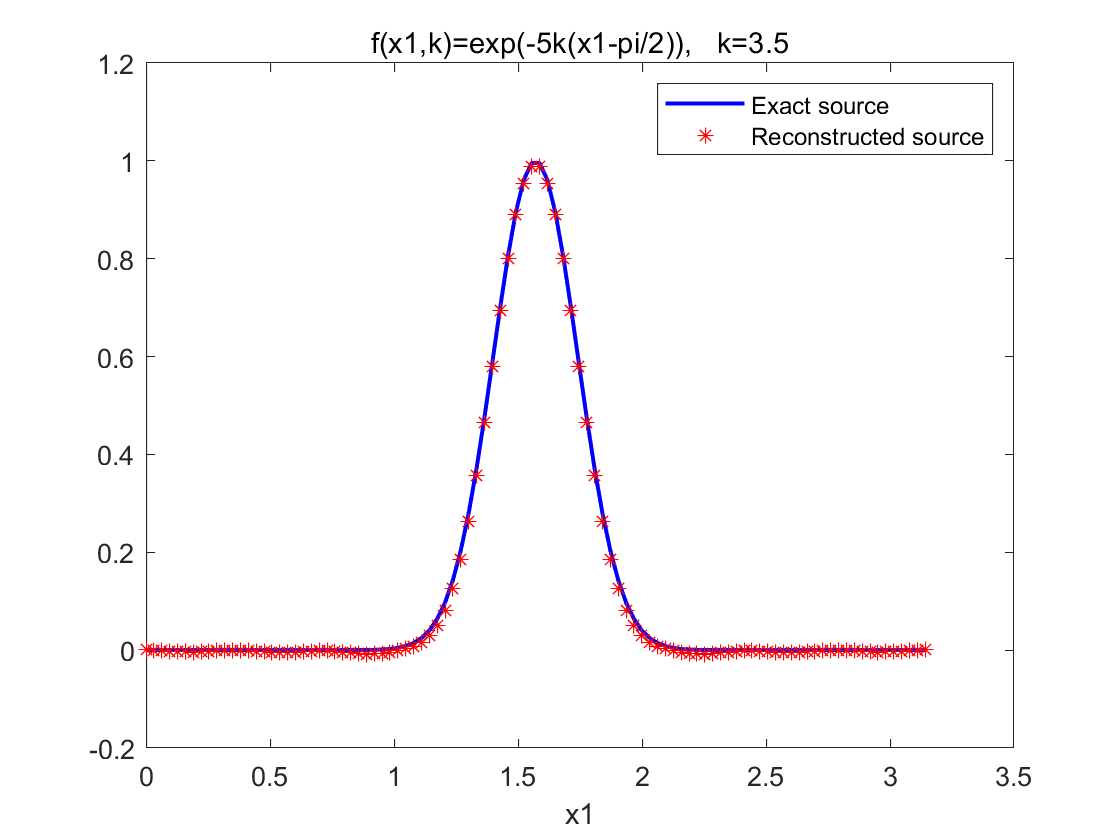}}
	\subfigure[$k=5.5$, $N=13$]{
		\includegraphics[scale=0.23]{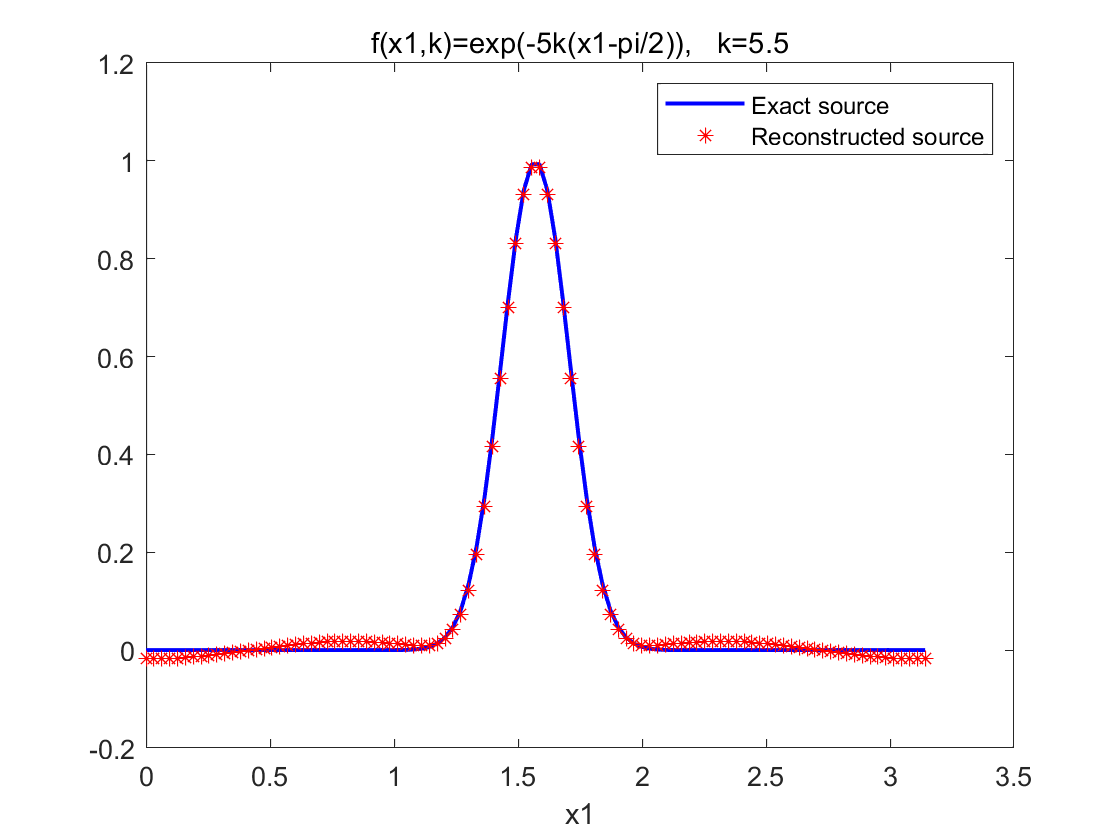}}
	\subfigure[$k=7.5$, $N=13$]{
		\includegraphics[scale=0.23]{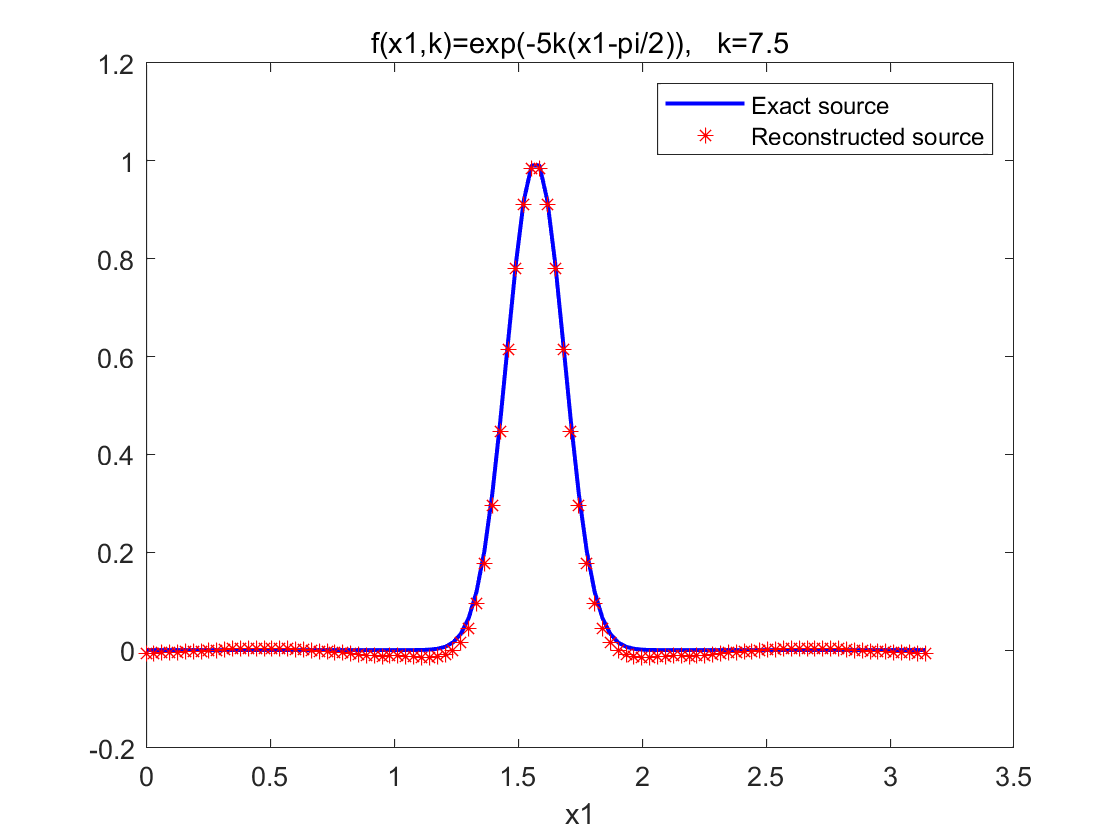}}		
	\caption{The reconstructed source using the Fourier-Transform method for the given source $e^{-5k(x_1-\pi/2)^2}$. Each figure corresponds to a different value of $k$.}
	\label{fig:FT:S1}
\end{figure}

It is important to note that for $k=2n$, $n\in \mathbb{N}^+$, the source function cannot be reconstructed due to the vanishing of $\tilde{G}_{n,j}$. However, by approximating such wave-numbers $k=2n$ with $k=2n\pm\varepsilon$, $0<\varepsilon\ll1$, $n\in \mathbb{N}^+$, Fig. \ref{fig:FT:k} (b) displays a successful reconstruction of $\cos(kx_1)$ supported in $[\pi/4,3\pi/4]$ with $k=2$. 
Following this, we apply the Fourier-Transform method to the source $e^{-5k(x_1-\pi/2)^2}$ by using different wave numbers. The recovery results are presented in Fig.\ref{fig:FT:S1}. The recreated image in Fig.\ref{fig:FT:S1}(a) closely resembles the exact source function with $k=1.5$ and $N=8$. In Fig.\ref{fig:FT:S1}(b), employing $k=3.5$ and $N=8$ for the source reconstruction yields favorable results. Subsequently, Fig.\ref{fig:FT:S1} (c) and Fig.\ref{fig:FT:S1} (d) display the reconstructed sources with $k=5.5$, $N=13$ and $k=7.5$, $N=13$, respectively. It is evident from these figures that as the wave number $k$ increases, the recovery effect becomes slightly less accurate.

\begin{figure}[H]
	\centering	
	\subfigure[$\delta=0.5\%$, $N=2$, ${\rm error=0.0965}$]{
		\includegraphics[scale=0.23]{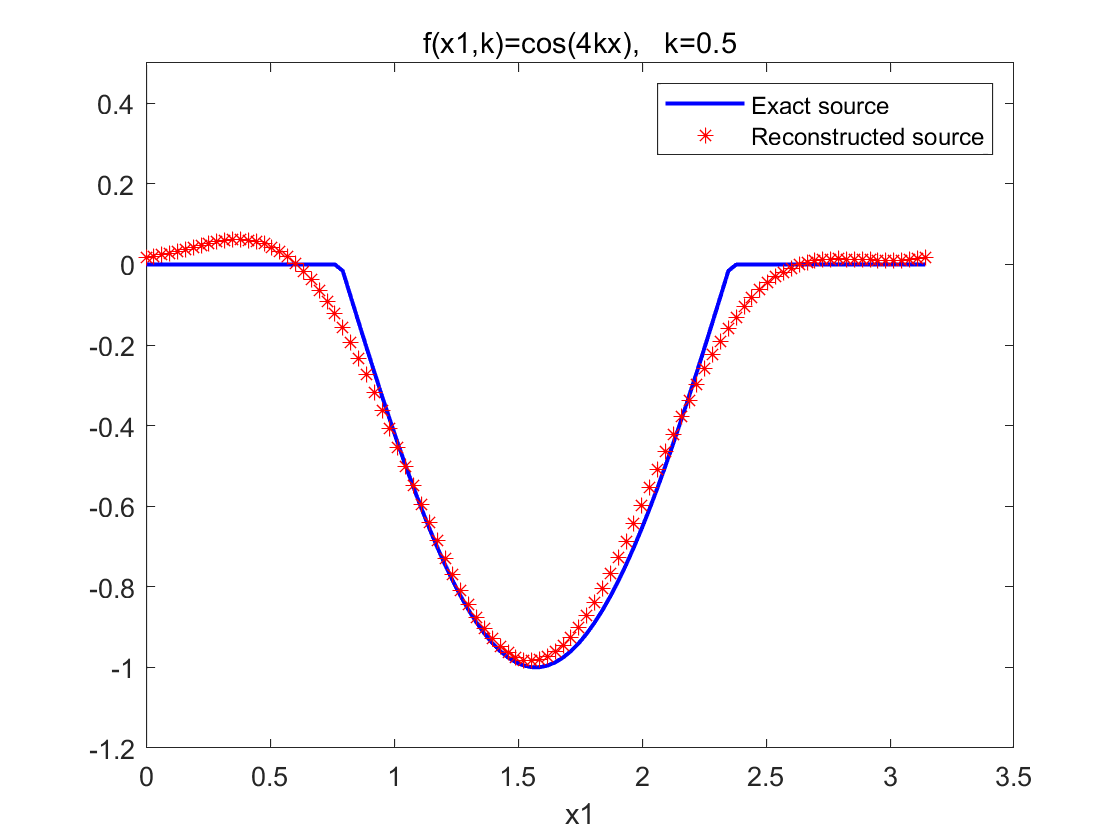}}
	\subfigure[$\delta=2\%$, $N=2$, ${\rm error=0.1637}$]{
		\includegraphics[scale=0.23]{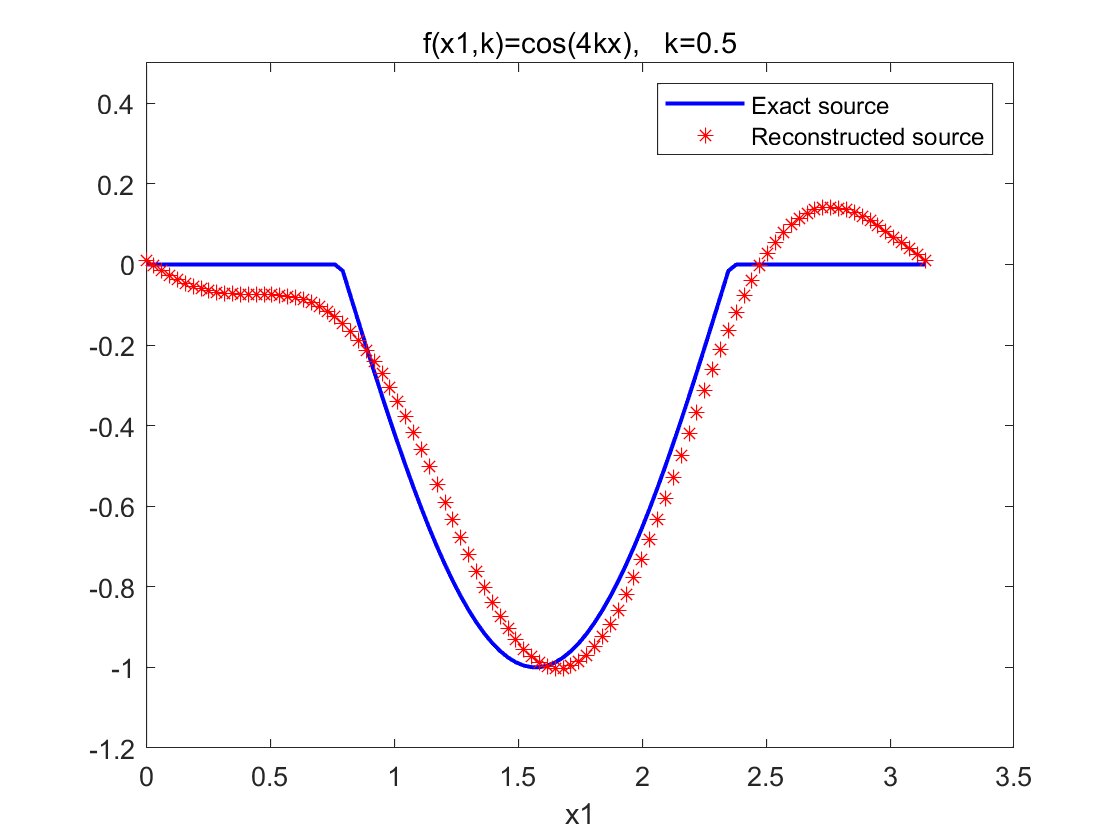}}	
	\subfigure[$\delta=10\%$, $N=1$, ${\rm error=0.3516}$]{
		\includegraphics[scale=0.23]{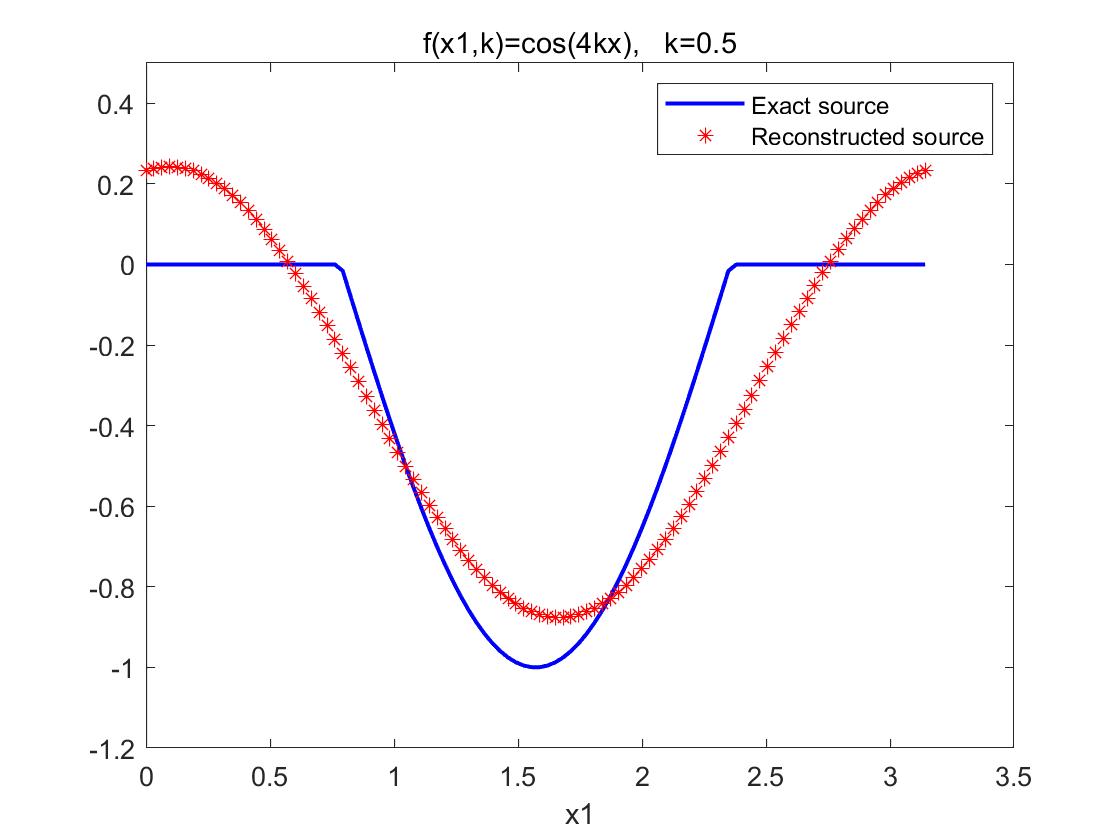}}
	\subfigure[$\delta=20\%$, $N=1$, ${\rm error=0.4106}$]{
		\includegraphics[scale=0.23]{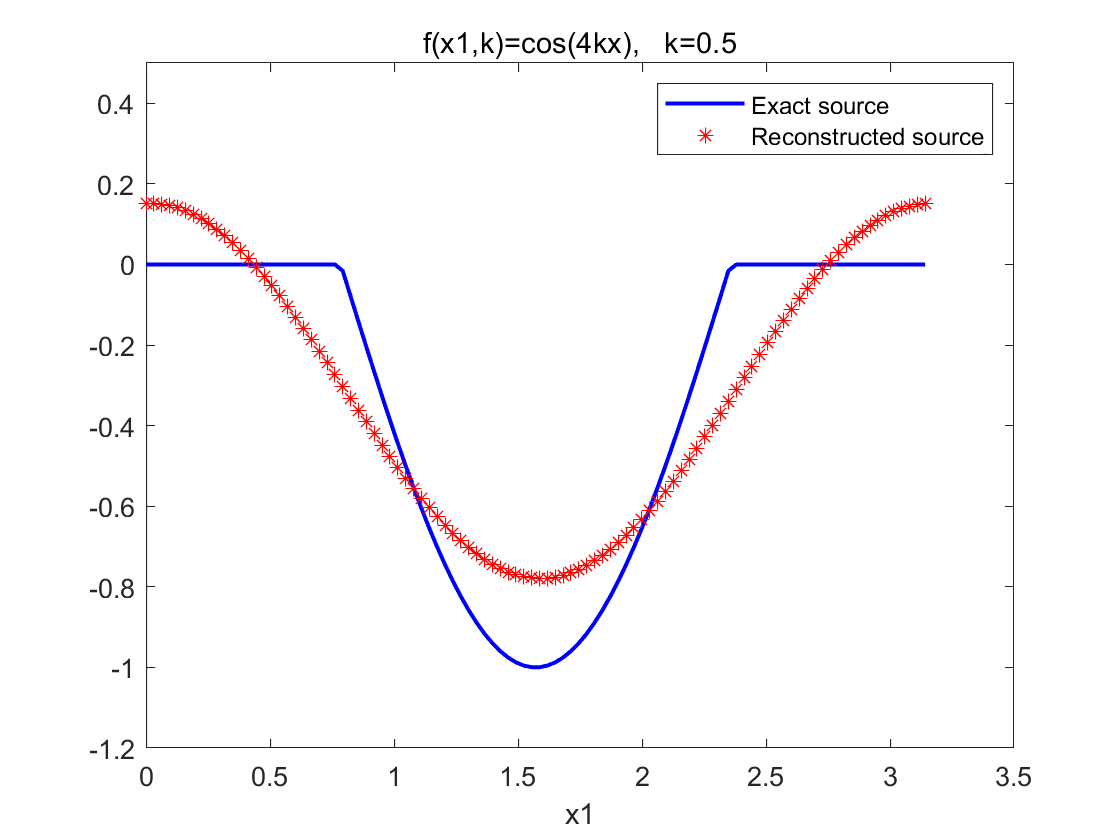}}
	
	\caption{The reconstructed sources using the Fourier-Transform method with different noise levels $\delta$.}
	\label{fig:FT:noise}
\end{figure}

At the end, we investigate the impact of noise on the Fourier-Transform method. The measurement data \eqref{equ:Ddata} is polluted by random noise using the formula \eqref{equ:noiseData}. Different noise levels, $\delta=0.5\%$, $\delta=2\%$, $\delta=10\%$ and $\delta=20\%$, are introduced. Fig.\ref{fig:FT:noise} illustrates the recreated results, indicating that the reconstruction error decreases as the noise level decreases, where the error is again defined by \eqref{equ:error}.

\section{Conclusion}

In this paper, we have established uniqueness, increasing stability and algorithms for an inverse wave-number-dependent source problem. In $d$-dimensions, the unknown source function depends on $\tilde{x}=(x_1,\cdots, x_{d-1})$ and $k$ but independent of $x_d$. The Fourier-Transform method, grounded in the Fourier basis, ensures both uniqueness and stability. Conversely, the Dirichlet-Laplacian method is constructed using Dirichlet-Laplacian eigenfunctions, though its stability remains to be futher examined. Two efficient non-iterative numerical algorithms have been developed. A possible continuation of this work is
to investigate  the stability under incomplete data and explore scenarios involving other kinds of wave-number-dependent sources.

\section*{Acknowledgements}
The work of G. Hu is partially supported by the National Natural Science Foundation of China
(No. 12071236) and the Fundamental Research Funds for Central Universities in China (No.
63233071). The work of S. Si is supported by the Natural Science Foundation of Shandong
Province, China (No. ZR202111240173).

\end{document}